\documentclass[12pt]{article}

\usepackage{latexsym,amssymb,amsfonts,euscript}
\usepackage{amsmath,amsthm,empheq,multirow}
\usepackage{graphicx}
\usepackage{subcaption}
\usepackage{adjustbox}
\usepackage{xcolor}
\usepackage[normalem]{ulem}
\usepackage{hyperref}

\newcommand{\bs}{\begin{eqnarray*}}
\newcommand{\es}{\end{eqnarray*}}

\newcommand{\field}[1]{\mathbb{#1}} % requires amsfonts
\newcommand{\bmat}{\begin{bmatrix} }
\newcommand{\emat}{\end{bmatrix}}
\newcommand{\beq}{\begin{eqnarray}}
\newcommand{\eeq}{\end{eqnarray}}
\newcommand{\beqs}{\begin{eqnarray*}}
\newcommand{\eeqs}{\end{eqnarray*}}

\newcommand{\goto}{\rightarrow}

\newcommand{\bm}{\mathbf}
\newtheorem{theorem}{Theorem}[section]

\definecolor{forestgreen}{rgb}{0.13, 0.55, 0.13}
\definecolor{darkmagenta}{rgb}{0.55, 0.0, 0.55}
\colorlet{colorYO}{darkmagenta}
\colorlet{colorGT}{blue}% color for Gamze
\colorlet{colorJM}{forestgreen}%color for JM

\newcommand{\JM}[1]{\textcolor{colorJM}{#1}}   % displays text in red color
\title{A Higher-Order Algorithm for Solving Time-Fractional Differential Equations via Rational Approximation of $z^{\alpha}$ as a Herglotz-Nevanlinna function}
\author{
Miao-jung Yvonne Ou \thanks{Department of Mathematical Sciences, University of Delaware, Newark, DE 19716, USA. \texttt{mou@udel.edu}, Corresponding author},
Gamze Tano\v{g}lu\thanks{Department of Mathematics, Izmir Institute of Technology, Izmir,T\"urkiye. \texttt{gamzetanoglu@iyte.edu.tr}},
Jiangming Xie \thanks{School of Mathematics and Statistics, Chongqing Technology and Business University, Chongqing, China. \texttt{xiejm@ctbu.edu.cn}}
}

\date{}

\begin{document}
\maketitle
\begin{abstract}
Time-Fractional derivatives (TFDs)  have been widely used in mathematical modeling of physical phenomena such as anomalous diffusion and memory effects of materials. However, an efficient time stepping scheme for solving Time-Fractional differential equations (TFDEs) remains a challenge due to the time-convolution term defining the TFD. Traditionally, TFD is handled {either by a \textbf{direct method} which saves almost all the memory and hence computationally expensive or a \textbf{sum of exponential method (SOE)} which decomposes the time domain improper integral into a local term and a history term} {approximated by} a sum of exponential functions. The limitation of the SOE approach is that { the approximation error of the local term might result in large errors if the unknown function is not smooth enough at $0$ and can lead to over-dissipation of energy.} In this paper, we take a different approach and will focus on the Caputo derivative 
because the Caputo time derivative of a constant function is 0. We present a  simple time-stepping algorithm based on the theory of Stieltjes-Markov functions, a subclass of the Herglotz-Nevanlinna functions. The key is in identifying the Stieltjes-Markov function structure hidden in the Laplace transform of the kernel function. Compared to most of the current methods in literature, the novelty of the results in this paper is manifold. First, no decomposition of integral is needed and our resulting approximation in the time domain contains a Dirac measure with carefully calculated strength, leading to a drastic decrease in the number of exponential functions needed in the approximation. Secondly, the choice of our interpolation points is based on the Mellin's inversion theorem and hence the interpolation error in the transformed domain is explicitly and rigorously linked to the error in the time-domain. {Finally, the resulting system is almost diagonal and hence can be easily sped up by using the Schur complement method.}

The performance of the proposed method is demonstrated through {four} test problems, including fractional ordinary differential equations with known analytical solutions, a fractional harmonic oscillator {and the nonlinear fractional Von der Pol equation}. The numerical experiments illustrate the accuracy and convergence properties of the method.
\end{abstract}

\section{Introduction}
Time-Fractional derivatives (TFDs)  have been studied extensively in the literature and 
     applied in modeling 
    {due to their ability of modeling innate memory effect in many physical phenomena such as tumor growth and wave propagation in viscoelastic media, see \cite{2010Fractional,metzler2000random,kilbas2006theory}.  A major theme in developing numerical schemes for TFDEs  has been tackling the difficulty of addressing
the nonlocal nature of the TFDs.   } 
The goal of this paper is to develop a fast method to handle the Caputo derivative of order $\alpha$, $0<\alpha<1$, defined in terms of the integral operator $\mathcal{I}^{1-\alpha}$ and the Gamma function $\Gamma(\cdot)$:
\beq
D^\alpha u(t):=\{\mathcal{I}^{1-\alpha}u' \}(t)&:=&\frac{1}{\Gamma(1-\alpha)} \int_0^t (t-s)^{-1+(1-\alpha)} u'(s)ds\nonumber\\
&=&\frac{1}{\Gamma(1-\alpha)} \int_0^t (t-s)^{-\alpha} u'(s)ds.
\label{Def:Caputo}
\eeq
{Further details concerning the definition is referred to \cite{diethelm2002analysis,podlubny1998fractional,samko1993fractional}.}

The Laplace transform of $D^\alpha u$, assuming $u(0)=0$, satisfies the following relation
\beq
\mathcal{L}[D^\alpha u](z)=\mathcal{L}\left[\frac{t^{-\alpha}}{\Gamma(1-\alpha)}\right](z) (\mathcal{L}[u](z)\cdot z -u(0))&=&\left(z^{\alpha-1}\right) \mathcal{L}[u](z)\cdot z  \nonumber\\
&=:& z^\beta \hat{u}\cdot z
\label{temp_Dalpha}
\eeq
where $-1<\beta:=\alpha-1<0$. {Note that $z^\beta$ corresponds to the Caputo fractional derivative $D^{1-\beta}$. }

{{Great efforts have been made to tackle the numerical challenge arising from the fractional derivative \eqref{Def:Caputo}. {{ In direct approaches \cite{jin2018discrete,lubich1988convolution,yuste2005explicit,zayernouri2016fractional,oldham1974fractional,zeng2015numerical}}}, the entire solution history is stored  and applied throughout the computation, which brings a heavy computation burden especially for long time simulations. Therefore, an efficient numerical method for handling the time convolution involved is of utmost importance. One of the popular techniques of handling the convolutional kernel is the summation of exponentials (SOEs), where the fractional derivative is replaced with local ordinary differential equations, {{e.g. \cite{Baffet2017,jiang2017fast,lubich2002fast,yuan2002numerical}}} There are also many other recent  works aiming for accelerating the convolution, {{\cite{banjai2019efficient,schadle2006fast,AghazadehMohammadiTanoglu2024,IdizTanogluAghazadeh2024,rosenfeld2017approximating,xie2019discontinuous,zheng2023fast} and references therein.  } } }

Despite the extensive study in the literature, an efficient numerical method for handling the time convolution involved is still an active research area. Most of the methods in the literature except \cite{10.1093/imanum/drac022,Baffet2017} have focused on approximating the kernel in the time domain with SOEs. In this paper, we take a different approach by first recognizing that the Laplace transform of the kernel function is a Herglotz-Nevanlinna function, which is expected because of the causal nature of the convolution term, followed by applying a simple linear fractional transformation to reveal the Stieltjes function structure of the Laplace-transformed kernel function. Using the well-known theorem of $[N-1/N]$ rational interpolation of Stieltjes functions \cite{gelfgren1978rational}, we interpolate the kernel function with grids located on a deformed Mellin's contour symmetric with respect to the real axis. This choice guarantees the poles to be all simple, negative with positive residues \cite{gelfgren1978rational}. The difference between our algorithm with what is  proposed in \cite{10.1093/imanum/drac022} will be explained in Section \ref{Conclusion}.

%{\color{blue}{In this paper, we take a different approach by first recognizing that the Laplace transform of the kernel function is a ....}}

To motivate the idea, we give a brief  description of the functions $z^\alpha$ with $-1<\alpha<1$ as Herglotz-Nevanlinna functions and their integral representation formulas (IRFs), which are crucial in building the numerical algorithm proposed in this paper.  According to p.218 in \cite{akhiezer1993theory-of-linea}, Theorem 1, A function $f$ defined and finite in the disc $|\zeta|<1$ admits the representation 
\[f(\zeta)=i\beta+\int_0^{2\pi} \frac{e^{it}+\zeta}{e^{it}-\zeta}\, d\sigma(t)
\]
where $\beta$ is a real constant and $\sigma(t)$ is a non-decreasing function (with bounded variation from the fact that $f(0)$ is finite) if and only if $f(\zeta)$ is regular and $Re(f)\ge 0$ in the disc $|\zeta|<1$, where $Re(\cdot)$ denotes the real part of a complex number.  

Then  with the transform
\[
z=i\frac{1+\zeta}{1-\zeta} \mbox{ and } \phi(z)=i f(\zeta)
\]
and rewriting the integral leads to the following characterization: A function $\phi(z)$ defined and finite in the half-plane $\mathcal{I}(z)>0$ admits the representation
\beq
\label{Herg_IRF}
\phi(z)=\alpha+\mu z+\int_{-\infty}^{\infty} \frac{1+tz}{t-z}\, d\sigma(t)
\eeq
where $\mu\ge 0$ and $\alpha$ are constants iff $\phi(z)$ is regular {(analytic)} in the upper half plane and has non-negative imaginary part. Such a function is called a Herglotz-Nevanlinna function and the following holds:
\beq
\label{two_const}
\alpha=Re(\phi(i)),\, \mu=\lim_{\eta\uparrow \infty}\frac{\phi(i\eta)}{i \eta}
\eeq

Another way to write \eqref{Herg_IRF} is
\beq
\phi(z)=\alpha+\mu z+\int_{-\infty}^{\infty} \left( \frac{1}{t-z} -\frac{t}{1+t^2}\right)\, d\tau(t)
\label{Nevanlinna_IRF}
\eeq
where the measures are related by $d\tau(t)=(1+t^2)d\sigma(t)$ and satisfies
\beq
\int_{-\infty}^\infty \frac{d\tau(t)}{1+t^2}<\infty
\label{measure_constraint}
\eeq

Now, consider the function with {the specific branch }
\beq
f(z)=z^{\alpha}:=|z|^{\alpha}e^{i\alpha\theta},\, -\pi<\theta:=Arg(z)\le \pi
\label{principal-branch}
\eeq
; this branch is chosen because we would like the poles and zeros to be all on the negative real axis. %for decaying in the time domain when Laplace transform is considered with $z=i\omega$.
Clearly, for $0<\alpha<1$, the imaginary part $\Im(z^{\alpha})=|z|^{\alpha}\sin(\alpha\theta)$ is positive if $0<\theta<\pi$ (upper half plane) and hence $z^\alpha$ and $-z^{-\alpha}$ are  Nevanlinna functions. Applying the Stieltjes inversion formula to \eqref{Nevanlinna_IRF}, the spectral function $\tau$ can be found explicitly. It is unique after the following normalization condition
\[
\tau(0)=0,\, \tau(t)=\frac{\tau(t+)+\tau(t-)}{2}
\]

For example, for $0<\alpha<1$, we have
\beq
z^\alpha&=&\cos\left(\frac{\alpha\pi
}{2}\right)+\frac{\sin(\alpha\pi)}{\pi}\int^0_{-\infty} dt (-t)^\alpha \left(\frac{1}{t-z} -\frac{t}{1+t^2}\right) \label{zS-1}\\
%&=& \cos\left(\frac{r\pi}{2}\right)+\frac{\sin{r\pi}}{\pi}\int^0_{-\infty} dt \frac{(-t)^r}{1+t^2} \left(\frac{1+tz}{t-z}\right)
z^{-\alpha}&=&\cos\left(\frac{\alpha\pi}{2} \right)-\frac{\sin(\alpha\pi)}{\pi}\int^0_{-\infty} dt (-t)^{-\alpha} \left(\frac{1}{t-z} -\frac{t}{1+t^2}\right) \label{zS}
\eeq

Furthermore, the integral representation formulas \eqref{zS-1}-\eqref{zS} can be brought into more compact forms \eqref{eq:znegdeltaOu} and \eqref{eq:zdeltaOu} by using the  theorems about the Stieltjes classes of function $\mathbf{S}$  and $\mathbf{S^{-1}}$ \cite{VKStieltjesFunctionsandHurwitz} and Lemma S1.5.2 in \cite{Kac1974R-functions-ana}. 

For $0<\alpha<1$, for $z\in \field{C}\setminus(-\infty,0]$, the following integral representation formulas (IRFs) are valid.
\begin{eqnarray}
\label{eq:znegdeltaOu}
z^{-\alpha}&=&\frac{\sin(\alpha\pi)}{\pi}\int_0^{\infty}  \lambda^{-\alpha} \left(\frac{1}{\lambda+z} \right)d\lambda,\\
z^{\alpha}&=&\frac{\sin((1-\alpha)\pi)}{\pi}\int_0^{\infty}  \lambda^{\alpha-1} \left(\frac{z}{\lambda+z} \right)d\lambda \label{eq:zdeltaOu}
\end{eqnarray}

What these IRFs reveal is that  rational function approximations are very natural for $z^{-\alpha}$ and $z^\alpha$ because the measures in the IRFs can be approximated by a finite sum of Dirac measures. However, the corresponding measure does not have finite mass and this prevent the use of many powerful tools tailored for subclasses of Herglotz-Nevanlinna functions whose IRFs correspond to measure with finite mass. One such subclass is the Stieltjes functions. In the next section, we will show that for $0<\alpha<1$, $z^\alpha$ and $z^{-\alpha}$ can be expressed in terms of a Stieltjes function and a constant. This observation is the foundation of the numerical algorithm proposed in Section \ref{sec:algorithm}.  

\section{The Stieltjes (Markov)  functions in the IRF}
{%\color{red}
We start with the characterization theorem of Stieltjes functions.

\begin{theorem}[Integral Representation of Stieltjes function \cite{akhiezer1993theory-of-linea}]
\label{Stieltjes-thm}
Let  $m(z)$ be a {Herglotz-Nevanlinna} function. Then  the following holds:

 \beq
 \sup_{y>0} y | m(iy) | <\infty \Longleftrightarrow
m(z)=\int_\field{R} \frac{d\mu(\lambda)}{\lambda-z}  \mbox{ and } \int_\field{R} d\mu(\lambda)<\infty 
\label{def:Stieltjes}
\eeq
A function satisfying the two integral conditions on the right of \eqref{def:Stieltjes} is called a Stieltjes function with its negative referred to as a Markov function. 
{Note: ${\overline{m(z)}=m(\overline{z})} $ }
\end{theorem}

%Partially due to  the fact that the IRF of this special class of Herglotz-Nevanlinna functions corresponds to measure of finite mass, 
The rational approximation of the Stieltjes functions is well studied \cite{gelfgren1981multipoint}. To utilize these powerful results, the focus of this section is to express \eqref{eq:znegdeltaOu} and \eqref{eq:zdeltaOu} in terms of Stieltjes functions through the following steps.. 
%With the choice of  branch-cut \eqref{principal-branch},  the function in \eqref{eq:znegdeltaOu} and \eqref{eq:zdeltaOu} have branch-cuts on the negative real axis, which is clearly unbounded. 
%On the other hand
%, to use the two-sided residue approximation described in Ou-Woederman, 
%we need the branch-cut to located on the positive real-axis.
%To map the branch cut to a finite interval on the positive real axis, we consider the new variable defined by the  following linear fractional transformation 
Consider the change of variable defined by 
\beq
s=\frac{1}{1-z}
\label{s_from_z}
\eeq

Under the action of  \eqref{s_from_z},  the upper half of the $z$-plane and  the branch-cut $z\in(-\infty,0]$ are mapped to the upper half of the $s$-plane and $s\in(0,1]$, respectively.
%%%%%%%%%%%%%%%%%%%%%%%%%%%%%%%%%%%%%%% 
\subsection{$z^\alpha$, $0<\alpha<1$}

Defining this way, the branch point of $z^\alpha$ at $z=\infty$  is mapped to $s=0$, 
\beq
z^\alpha=\left(1-\frac{1}{s} \right)^\alpha=:G(s)
\label{Gs}
\eeq
%We should try to find the IRF of $G(s)$, what will be the corresponding values of $a$ and $b$? Will this give us a non-zero $a$?
%%%%%%%%%%%%%
\begin{theorem}
There exists a positive Borel measure $d\sigma(t)$ with finite mass such that the $G(s)$ in \eqref{Gs}  can be expressed as follows 
\beq
G(s)&=&1+ \int_0^1  \frac{d\sigma(t)}{t-s}, \mbox{ for } s \in \field{C}\setminus[0,1]\\
&\approx&1+ \sum_{k=1}^N \frac{r_k}{p_k-s},\quad r_k>0,\quad p_k\in(0,1),
\label{approx_P}
\eeq
where $p_k\ne p_j$ if $k\ne j$.
\end{theorem}
\begin{proof}
Note that $G(s)$ approaches $1$ as $s\goto \infty$. 
%%%%%%%%%%%%%%%%%%%%%%%%%%%%%%%%%%%%%%%%%%%%%%%%%

%%%%%%%%%%%%%%%%%%%%%%%%%%%%%%%%%%%%%%%%%%%%%%%%%%%%%%%%%
With a simple geometric argument, we see that  $G(s)$ maps $C_+$ to $C_+$ under the selected branch  inherited from \eqref{principal-branch} for $0<\alpha<1$.

Define the new function 
\beq
P(s):=G(s)-1=\left(1-\frac{1}{s} \right)^\alpha-1%,\quad s=\frac{1}{1-z},\mbox{ with the selected branches.}
\label{def_P}
\eeq 
%Then for $s\in \field{C}_+$, $P(s)$ has the following Taylor series expansion
%\beq
%P(s+\delta)=P(s)+P'(s)\delta+O(\delta^2)=P(s)+\alpha s^{-2} (1+\frac{1}{s})^{\alpha-1}\delta+O(\delta^2). 
%\eeq
It is important to note that for $s=iy, y>0$, the following holds:
\beqs
\lim_{\field{R}\ni y \goto \infty} yP(iy)=\lim_{\field{R}\ni y \goto \infty} \frac{(1+ i \frac{1}{y})^\alpha -1}{\frac{1}{y}}=\lim_{y \goto 0^+} \frac{(1+ i y)^\alpha -1}{{y}}.
\eeqs
Since the function $f(\xi):=(1+\xi)^\alpha$ is analytic for $\xi\in\field{C}\setminus \{\xi\le -1\}$, the Taylor expansion of $f(\xi)$  at $\xi=0$ leads to 
\beqs
\lim_{\field{R}\ni y \goto \infty} yP(iy)=
\lim_{y \goto 0^+} \frac{1+i\alpha y +O(y^2) -1}{{y}}=i\alpha.
\eeqs
Moreover, 
\beqs
\lim_{y \goto 0^+} yP(iy) =  \lim_{\field{R}\ni y \goto \infty} \frac{\left( \sqrt{1+y^2}\right)^\alpha\cdot e^{i\alpha {\rm Arg}(1+iy)}-1 }{y}=0, \mbox{ because $0<\alpha<1$.}
\eeqs
Hence for $0<\alpha<1$, $|yP(iy)|$ is bounded for all $y>0$ and therefore Theorem \ref{Stieltjes-thm} implies that $P(s)$ is a Stieltjes function analytic in $\field{C}\setminus [0,1]$. The existence of approximation \eqref{approx_P} with the special properties of the poles $p_k$ and the residues $r_k$ is a direct application of the theorems of Markov functions in \cite{gelfgren1978rational}.
\end{proof}

Equation \eqref{approx_P} leads to the following estimation:
\beq
z^\alpha \approx \left(1+\sum_k \frac{r_k}{p_k} \right) - \sum_k \frac{r_k/p_k^2}{z-( 1-\frac{1}{p_k})}
\label{approx_z_pos}
\eeq
Note that $1-\frac{1}{p_k}<0$ because $p_k\in (0,1]$. 
\subsection{$z^{-\alpha}$, $0<\alpha<1$}
In terms of variable $s$ defined in  \eqref{s_from_z}, we have  
\beq
-z^{-\alpha}=-\left(1-\frac{1}{s} \right)^{-\alpha}=: H(s), 
\eeq
%%%
\begin{theorem}
\label{Thm_Q}
There exists a positive Borel measure $d\mu(t)$ of finite mass such that
\beq
H(s)=\int_0^1 \frac{d\mu(t)}{t-s}-1 \approx\sum_{k=1}^N \frac{ {r}_k}{ {p}_k-s}-1,\,  {r}_k>0,\, {p}_k\in (0,1)
\label{approx_Q}
\eeq
where $p_k\ne p_j$ if $k\ne j$.
\end{theorem}
\begin{proof}
Note that $H(s)$ maps from $\field{C}_+$ to $\field{C}_+$. Since $H\goto -1$ as $s\goto \infty$, we define the new function 
\beq
Q(s):=H(s)+1=-\left(1-\frac{1}{s} \right)^{-\alpha}+1.
\label{def_Q}
\eeq 
Then it is easy to show that  $Q(s)$ maps $\field{C}_+$ to $\field{C}_+$, approaches 0 as $s\goto 0$, and is analytic for $s\in\field{C}\setminus (0,1]$. Moreover, we have for $y>0$,
\beqs
\lim_{\field{R}\ni y \goto \infty} yQ(iy)=\lim_{\field{R}\ni y \goto \infty}  \frac{-(1+ i \frac{1}{y})^{-\alpha} +1}{\frac{1}{y}}=\lim_{y \goto 0^+}\frac{-(1+ i y)^{-\alpha} +1}{{y}}=i\alpha 
\eeqs
Also,

\beqs
\lim_{y \goto 0^+} yQ(iy) = \lim_{0<y \goto \infty} \frac{[ -(1+ i y)^{-\alpha} +1]}{y}=\lim_{0<y \goto \infty} \frac{-\left( \sqrt{1+y^2}\right)^{-\alpha}\cdot e^{-i\alpha {\rm Arg}(1+iy)}+1 }{y}=0.
\eeqs
Hence for $0<\alpha<1$, $|yQ(iy)|$ is bounded for all $y>0$ and therefore $Q(s)$ is a Stieltjes function by Theorem \ref{Stieltjes-thm} and hence the IRF in \eqref{approx_Q} is valid. Moreover, the existence of approximation \eqref{approx_Q} with the special properties of the poles $p_k$ and the residues $r_k$ is a direct application of the theorems of Markov functions in \cite{gelfgren1978rational}. %
\end{proof}

Theorem \ref{Thm_Q} guarantees the existence of the rational approximation specified below 
\beq
z^{-\alpha}&=&1-Q(s) \approx \left(1-\sum_{k=1}^N \frac{ {r}_k}{ {p}_k}\right)+\sum_{k=1}^N \frac{ {r}_k/ {p}_k^2}{z-\left(1-\frac{1}{ {p}_k}\right)}\nonumber\\
&=:&  {c}_N +\sum_{k=1}^N \frac{ {\gamma}_k}{z- {\eta}_k}=:R_{N,-\alpha}(z),\,  {\gamma}_k>0,\,  {\eta}_k<0,
\label{approx_z_neg}
\eeq
where $ {c}_N:=1-\sum_{k=1}^N \frac{ {r}_k}{ {p}_k}$, $ {\gamma}_k:= {r}_k/ {p}_k^2$ and $ {\eta}_k:=1-\frac{1}{ {p}_k}$. Note that $ {p}_k\in (0,1]$ so $ {\eta}_k<0$. Therefore, in the time domain, this results in \eqref{approx_time_domain} consisting of an approximates the time convolution kernel with an instantaneous term $c_N u'(t)$, a dissipation term $D_N u(t)$ and $N$ relaxation terms with distinct relaxation times.

%We would like to remark that the IRF of $P(s)$ and $Q(s)$ provide the theoretical foundation for the existence of rational approximations with the forms specified in \eqref{approx_z_pos} and \eqref{approx_z_neg}. That is, all the poles must be simple and locate on the negative real axis with positive residue; it can be realized by applying the theorems from \cite{gelfgren1978rational} with  properly chosen  grid points. 
Contrary to what is proposed in\cite{10.1093/imanum/drac022}, it is important to note that the constants in  \eqref{approx_z_pos}  and \eqref{approx_z_neg} depend on $\{r_k, p_k\}_{k=1}^N$ and hence should not be treated as independent parameters if optimization schemes are used to find the coefficients of the rational approximation.

%%%%%%%%%%%%%%%%%%%%%%%%%%
\section{Rational approximations through interpolation of the Stieltjes functions}
The approximations in \eqref{approx_z_pos} and \eqref{approx_z_neg} rely on the approximations of the Stieltjes functions in \eqref{approx_P} and \eqref{approx_Q}, respectively, and we refer to them as the 'pole-residue' approximations. To find the poles and residues such that all poles are simple and located in $(0,1)$ with positive residues, we apply the interpolation results of Markov functions in \cite{gelfgren1978rational}. Before describing the algorithms and lemmas for constructing the poles and residues, we would like to remark on the approximation errors between the $z$-domain and the time domain $t$. 
\subsection{\label{error-link}The link between the approximation errors in the $z$-domain and  the time-domain}
Consider the pointwise error of the pole-residue approximation of $z^{-\alpha} $ 
\[
E_N(z):=z^{-\alpha}-R_{N,-\alpha}(z)=Q_N(s)-Q(s)
\]
where $Q(s)$ is defined in \eqref{def_Q} and $Q_N(s):=\sum_{k=1}^N \frac{ {r}_k}{ {p}_k-s}$ is a rational approximation of $Q(s)$.
{
A direct application of Theorem 5.1 in \cite{gelfgren1978rational} gives the  interpolation error estimate for $z\in \mathcal{C}_L$ as follows. 
\begin{theorem}
\label{thm_N_conv}
%If $a,b$ are finite numbers, for $s\in \field{C}\setminus[a,b]$ with $s_k$, $k=1,\cdots 2n$, as the interpolation points, 
\[ 
\left|E_N(z)\right| \le K_v\cdot \frac{1}{max(|s|,|s-1|)}\cdot \frac{1}{(1-|\psi_v(s)|)} \cdot \Pi_{k=1}^{2N} |G_k(s)|
\]
where $s=\frac{1}{1-z}$, $\psi_v:=\frac{     \sqrt{s-1} -  \sqrt{\frac{v-1}{v}}\cdot \sqrt{s}   }{     \sqrt{s-1} +  \sqrt{\frac{v-1}{v}}\cdot \sqrt{s}   }$,  a conformal mapping of $\field{C}\setminus[0,1]$ onto the interior of the unit circle, $v=s_j$ for some $j$ and  $G_k(s):=\frac{\psi_v(s_k)-\psi_v(s)}{1-\overline{\psi_v(s_k)} \psi_v(s)}$,  $K_v$ depends on $v$ and $\alpha$ but not on ${\JM{N}}$.
\end{theorem}
Note that $|G_k(s)|<1$ for all $z\in\mathbb{C}\setminus[0,1]$ because of the mapping property of $\psi_v$ and the maximum modulus principle for functions analytic inside the unit disk; hence $|E_N(z)|$ converges to $0$ geometrically as $N$ increases.
}

The ultimate question is: "How does the approximation error on the $z$-plane (or the $s$-plane) propagate to the time domain?", which can be addressed by the following theorem.

%Mellin's inversion formula for the Laplace transform 

\begin{theorem}[Mellin's inversion formula  (Exercise 11.39 in \cite{apostol1974mathematical-an})] 
\beq
\frac{f(t+)+f(t-)}{2}=\frac{1}{2\pi i}\lim_{y\goto \infty} \int_{\eta-iy}^{\eta+iy} e^{zt} F(z) dz,\, \eta>c
\label{Mellin}
\eeq
if $f$ satisfies at least one of the two local conditions:
\begin{itemize}
\item
There exist $\delta>0$ such that $f$ is of bounded variation on $[t-\delta,t+\delta]$
\item
Both limits $f(t+)$ and $f(t-)$ exist and both of the Lebesgue integrals $\int_0^\delta \frac{f(t+\tau)-f(t+)}{\tau}\, d\tau$ and $\int_0^\delta \frac{f(t-\tau)-f(t-)}{\tau}\, d\tau$ exist
\end{itemize}
and its Laplace transform $F(z)$ exists for $Re(z)>\eta$. 
\end{theorem}
In this paper, we will refer the contour used in \eqref{Mellin} as the Mellin contour; see Figure \ref{Parabola}.

Because of  the constant $c_N$ in \eqref{approx_z_neg},  the inverse transform of $E(z)$ contains a Dirac measure $c_N \delta(t)$ and hence cannot be compared with other functions in a point-wise sense. Therefore, we will treat both $z^{-\alpha}$ and its approximation as a multiplication operator acting on a function $U(z)$ whose singularities are all to the left of a vertical line $Re(z)=\eta$ for some $\eta>0$. This viewpoint is motivated by the fact that the fractional derivative corresponds to a product in the $z$-domain, and a good approximation along a curve of the kernel function  alone does not necessarily lead to a good approximation of the convolution unless the inverse Laplace transform of the product evaluated along the curve really corresponds to the convolution integral of $D^{1-\alpha} u$, i.e. all the singularities of the function $U(z)$ are also to the left of the curve. 

Noting that if $U(z)$ has all its singularities located to the left of a parabola $x=\eta-y^2$ shown in Figure \ref{Parabola}, then $E(z)$ is analytic in $\field{C}\setminus\{z<0\}$ and hence the Mellin contour can be replaced by this parabola according to the Cauchy integral theorem. \footnote{This was suggested to us by Prof. Andr\'e Weideman (Stellenbosch University, South Africa, during the EXPOWER workshop in University of Lille in summer 2024.} 
%This choice of contour is to address the fact that we can only interpolate on a finite portion of the parabola and have no control of the error $E(z)$ outside it  so we would like the contribution of error from this part of the contour to damp to 0 fast with the help of the factor $e^{zt}$. 
This type of contours (Hankel contour) has been used in the numerical Laplace transform literature because it results in a fast convergence of the integral in \eqref{Mellin}; see \cite{dingfelder2015improved} and the references therein.

As will be seen in the analysis below, for $t>0$, this choice of contour will make the error outside the interpolation region very small because of the Mellin kernel $e^{zt}$. Introduce the notation of truncated parabola 
\beq
\mathcal{C}_L:=\{z=(\eta-y^2)+iy,\quad -L\le y\le L \}
\label{eq:CL}
\eeq
; see Figure \ref{Parabola}. 
%%%%%%%%%%%%%%%%%%%%%%%%%%
%%%  This figure is added by JM%%%
\begin{figure}[!h]
\begin{center}
\includegraphics[width=0.5\textwidth, height=0.3\textheight]{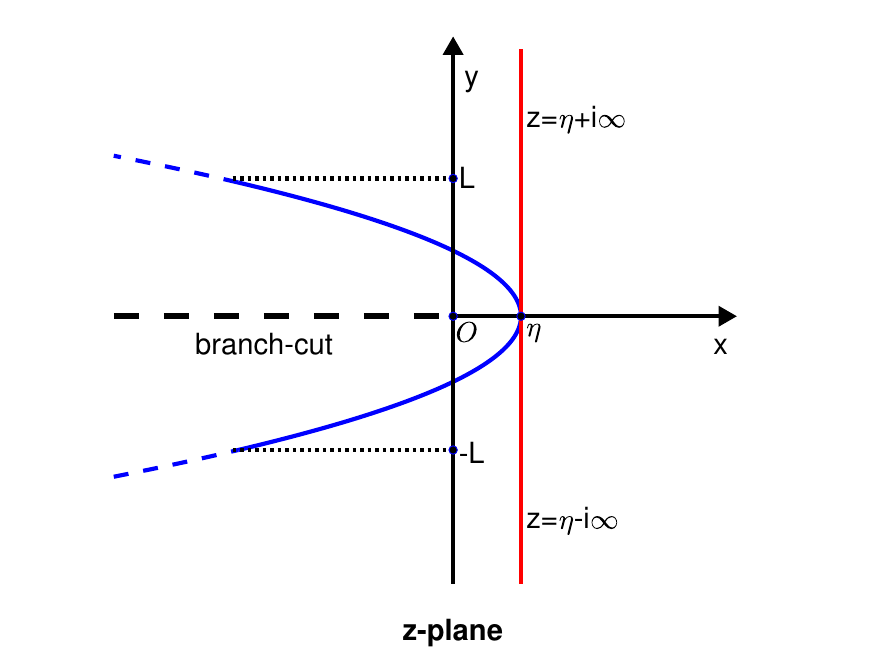}
\caption{Red: the Mellin contour. Blue: The deformed contour. Blue solid curve: $\mathcal{C}_L$ \label{Parabola}}
\end{center}
\end{figure}
%%%%%%%%%%%%%%%%%%%%%%%%%
%
%%%%%%%%%%%%%%%%%%%%%%%%%%
The interpolation grids will be on $\mathcal{C}_Y$ with some pre-defined constant $Y>0$. To utilize the theorems of Stieltjes functions, the interpolation will be carried out in the $s$ domain in order to first construct $Q_N(s)$ with the corresponding grids   via \eqref{s_from_z}. To be more precise,  $Q_N$ is a rational function of type $[N-1/N]=\frac{p_{N-1}(s)}{q_N(s)}$, where $p_{N-1}(s)$ and $q_N(s)$ are polynomials of degree $N-1$ and $N$, respectively. So there will be a total of $2N$ interpolation points in $\field{C}\setminus[0,1]$ on the $s$-plane. In order to guarantee that all poles are simple and located in $(0,1)$ with positive residue, the complex-valued interpolation points need to come in complex-conjugate pairs; we denote these $k_2$ pairs of grids by $\{z_j, \overline{z_j}\}_{j=1}^{k_2}$. With the real-valued grid $x_1:=\eta$, we need to also interpolate at one additional point to make the total number of grids an even number. Following Theorem 1 in \cite{gelfgren1978rational}, we choose this addition point to be $s=\infty$. 

Specifically, the interpolation points on the $s$-plane are selected as follows:
\beq
&&s_j=\frac{1}{1-(\eta-y_j^2+iy_j)},\quad j=1,\dots,N-1=:k_2 ,\nonumber\\
&& y_j=j\left(\frac{Y}{N-1}\right) \mbox{ for equal spacing grid }; y_j=e^{j\left(\frac{\ln(Y)}{N-1}\right)} \mbox{ for log spacing grid }\nonumber\\
&&s_{j+k_2}=\overline{s_j}, \quad j=1,\dots,k_2,\nonumber\\
&&s_{2N-1}=\frac{1}{1-\eta}, \quad s_{2N}=\infty,
\label{grids}
\eeq

We would like to remark that the popular two-sided residue approximation (TSRA)\cite{Daniel-Alpay1994The-two-sided-r} is not applicable for this problem because TSRA only handles grids which are  complex-conjugate  pairs. 

{
The time domain error corresponding to  $R_{N,-\alpha}$ is
\beqs
\mathcal{E}(t;u)&:=& \frac{1}{2\pi i}\lim_{L\rightarrow \infty}\int_{\mathcal{C}_L} (R_N(z)-z^{-\alpha})\, U(z)e^{tz}\, dz %\equiv \frac{1}{2\pi i} \int_{\mathcal{C}_\infty} (R_N(z)-z^{-\alpha})\, U(z)e^{tz}\, dz 
\eeqs
such that}
{
\beq
|\mathcal{E}(t;u) |&\le& \frac{1}{2\pi} \left[\int_{\mathcal{C}_\infty} |R_{N,-\alpha}-z^{-\alpha}|^2 |e^{tz}| |dz|\right]^{\frac{1}{2}}
\left[\int_{\mathcal{C}_\infty} |U(z)|^2 |e^{tz}| |dz|\right]^{\frac{1}{2}}\nonumber\\
&=& \frac{e^{\eta t}}{\pi}\left[ \int_0^Y \epsilon^2(y)e^{-ty^2} ds(y)  \right. \nonumber\\
&&+ \left. \int_Y^\infty |R_{N,-\alpha}(z(y))-z^{-\alpha}(y)|^2 e^{-ty^2} ds(y)\right]^{\frac{1}{2}} 
\cdot \| U\|_{\mathcal{C}_\infty}(t)
%&=:& \mathcal{E}_{|y|>Y}(t;u)+\mathcal{E}_Y(t;u),\quad t>0,
\label{error_estimate}
\eeq
where $ds(y)=\sqrt{1+4y^2}\,dy$, $\epsilon(y)=|E_N(z(y))|$ and 
}
{
\[
\| U\|_{\mathcal{C}_\infty}(t)=\left[\int_0^\infty |U(z(y))|^2 e^{-ty^2} ds(y)\right]^{1/2}
\]
Note that the first term in the brackets in \eqref{error_estimate} corresponds to the interpolation error while the second term the extrapolation error all the terms there decay exponentially as time $t$ increases except the factor $e^{\eta t}$. In practice, if we are to solve the TFDs up to $t=T$, we can set $\eta=\frac{1}{T}$. Also, the analysis in \eqref{error_estimate} shows the clear advantage of using the deformed contour over the standard Mellin contour for interpolation of the kernel. }

} }
{Based on the error formula \eqref{error_estimate}, we see that the error of the approximation is larger for very small $t$, and this inspire the shifted scheme; see \cite{Baffet2017} and the references therein.}%
\begin{figure}[htbp]
\centering
    \begin{subfigure}{0.48\linewidth}
        \centering
        \includegraphics[width=\linewidth]{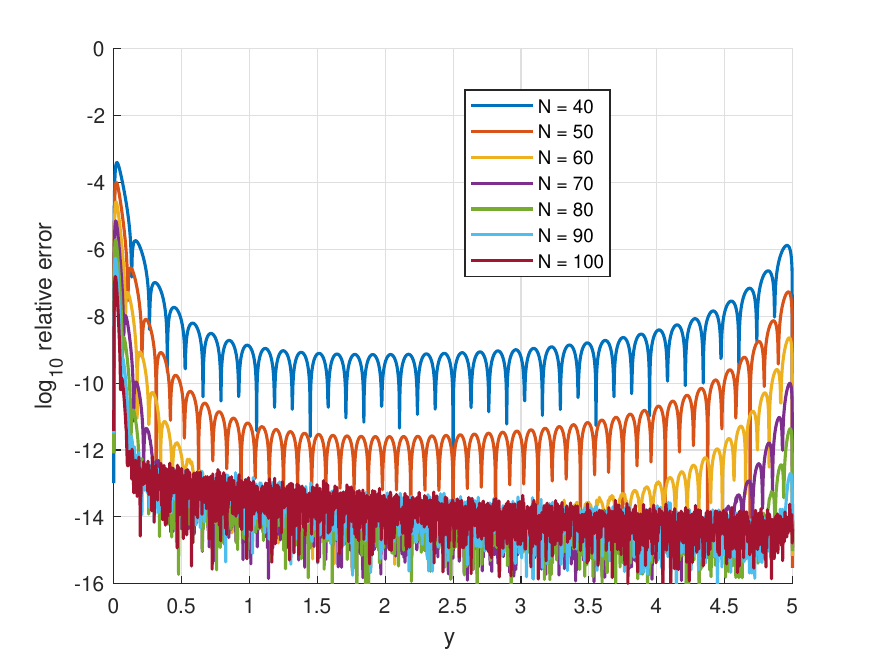}
        \caption{}
    \end{subfigure}\hfill
    \begin{subfigure}{0.48\linewidth}
        \centering
        \includegraphics[width=\linewidth]{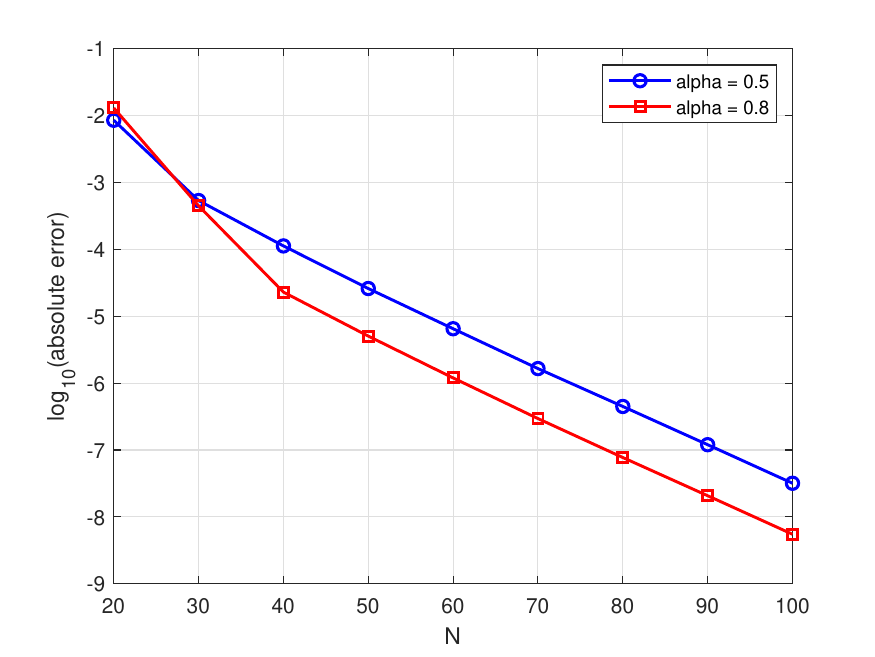}
        \caption{}
    \end{subfigure}
    \caption{Numerical errors of the rational approximation of $z^{-\alpha}$. (a) the interpolation error $E_N$ with $\alpha=0.8$ and  $N=40,\,50,\,...,\,100$, (b) the maximum error of $E_N$ with different $\alpha$ and $N$. \label{fig:numericalE_N}}
\end{figure}
%%%%%%%%%%%%%%%%%%%%%%%%%%
%%%%%%%%%%%%%%%%%%%%%%%%%%
\begin{figure}[!b]
\begin{center}
\includegraphics[width=0.9\textwidth, height=0.5\textheight]{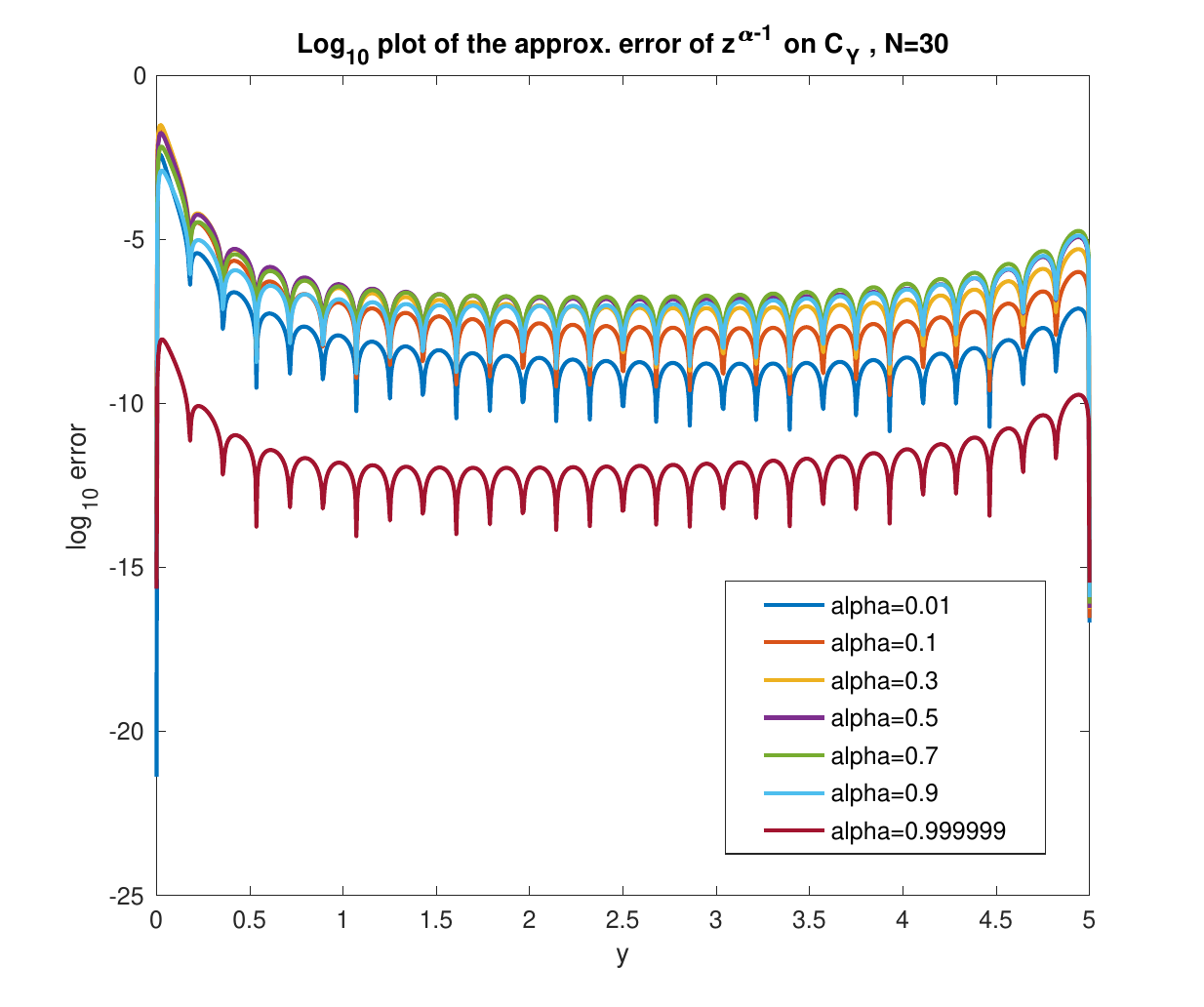}
\caption{Log-10 of the interpolation error $\epsilon(y)$ for various values of $\alpha$, where $\eta=0.01$ and $N=30$  \label{compare}}
\end{center}
\end{figure}
%%%%%%%%%%%%%%%%%%%%%%%%%%

%maximum error of the approximation of $z^{-\alpha}$ by $F$ on the truncated parabola $\mathcal{C}_Y$, and ${\rm erf}(\cdot)$ the error function. Usually, we observe the error dropped quickly in a boundary layer near $y=0$, therefore, we will let $\epsilon_1$ denote the max. error in this boundary layer and $\epsilon_2$ for the rest; in general, $\epsilon_1>> \epsilon_2$, and hence a tighter estimate for $E_2$ can be obtained.

%For a predefined simulated time interval $[\delta,T]$, and the number of interpolation points $N$, \eqref{E1} and \eqref{E2} will lead to a point-wise error estimate for the convolution kernel and hence the error estimate on the solutions of the equations that contain a fractional time derivative term. 
%
\subsection{Linear system approach for constructing $Q_N$}
This method involves a linear system for solving the coefficients of the rational functions. It is straight forward but the linear system to be solved is very ill-conditioned. This approach consists of the following steps.
\begin{enumerate}
\item
Given the $2N$ grids $\{s_k\}_{k=1}^{2N}$ specified in \eqref{grids}. Define
\[
[N-1/N](s):=\frac{a_0+a_1s+\dots+ a_{N-1} s^{N-1}}{1+b_1 s+\dots+b_{N}s^{N}}
\]
\item
Evaluate the Stieltjes function $Q$ at $s_k$, $k=1,\dots, 2N-1$ 
\[
Q_k:=Q(s_k)=-\left(1-\frac{1}{s_k} \right)^{-\alpha}+1, \mbox{ with the  branch  in \eqref{principal-branch}}
\]
\item We note that
\beq
\lim_{s\goto \infty} sQ(s)=-\alpha.
\label{at_inf}
\eeq
To interpolate at $s_{2N}=\infty$ means
\beq
\lim_{s\goto \infty} s [N-1/N](s)=-\alpha
\eeq
and leads to the equation
\beq
a_{N-1}=-\alpha b_N
\label{last_eqn}
\eeq
\item
Solve the system of the following $2N-1$ equations
\beq
\frac{a_0+a_1s+\dots+ a_{N-1} s^{N-1}}{1+b_1 s+\dots+b_{N}s^{N}}=Q(s),\, s=\{s_k\}_{k=1}^{2N-1}, 
\eeq
and equation \eqref{last_eqn}.

The corresponding linear system is
\beq
\begin{pmatrix}
\mathbf{C}\\
\overline{\mathbf{C}}\\
\mathbf{u}\\
\mathbf{v}
\end{pmatrix}
\mathbf{x}=\begin{pmatrix}
\mathbf{d}\\
\overline{\mathbf{d}}\\
Q_{2N-1}\\0
\end{pmatrix}
\label{complex_system}
\eeq
where the overline represents the complex conjugate, and $\mathbf{C}$, $\mathbf{u}$, $\mathbf{v}$ and $\mathbf{d}$ are defined as follows:
\beqs
&&\mathbf{C}:=\begin{pmatrix}
1&s_1& s_1^2&\cdots s_1^{N-1} & -Q_1s_1 & -Q_1 s_1^2&\cdots& -Q_1s_1^{N}\\1&
s_2& s_2^2&\cdots s_2^{N-1} & -Q_2s_2 & -Q_2 s_2^2&\cdots& -Q_2s_2^{N}\\
\vdots& \vdots& \vdots&  \vdots& \vdots& \vdots& \vdots\\
1&s_{N-1}& s_{N-1}^2&\cdots s_{N-1}^{N-1} & -Q_{N-1}s_{N-1} & -Q_{N-1} s_{N-1}^2&\cdots& -Q_{N-1}s_{N-1}^{N}
\end{pmatrix}\\
&&\mathbf{u}=\begin{pmatrix}
\mathbf{u}_1 & \mathbf{u}_2
\end{pmatrix},\,
\mathbf{v}=\begin{pmatrix}
\mathbf{v}_1 & \mathbf{v}_2
\end{pmatrix},\,
\mathbf{x}=\begin{pmatrix}
a_0&a_1&\cdots&a_{N-1}&b_1&\cdots&b_{N}
\end{pmatrix}
\eeqs
where
\beqs
&&\mathbf{u}_1=\begin{pmatrix} 1& s_{2N-1}& s_{2N-1}^2&\cdots s_{2N-1}^{N-1}
\end{pmatrix}\\
&&\mathbf{u}_2=\begin{pmatrix}
-Q_{2N-1}s_{2N-1} & -Q_{2N-1} s_{2N-1}^2&\cdots& -Q_{2N-1} s_{2N-1}^{N}
\end{pmatrix}\\
&&\mathbf{v}_1=\begin{pmatrix} 0& 0& 0&\cdots &1
\end{pmatrix}\\
&&\mathbf{v}_2=\begin{pmatrix}
0 & 0& 0&\cdots& \alpha
\end{pmatrix}\\
&&\mathbf{d}=\begin{pmatrix}
Q_1& Q_2&\cdots& Q_{N-1}
    \end{pmatrix}^t
\eeqs
\item
According to Gelfgren's theorem \cite{gelfgren1978rational}, the interpolation grids chosen this way guarantee that the denominator of $[N-1/N]$ has $N$ distinct real-valued roots located in $(0,1)$. So we will apply partial fraction decomposition to find these roots $p_k$ and the corresponding weights $r_k$, $k=1,\dots, N$. This implies that the unknown coefficients are all real-valued. As is pointed out in \cite{ou2014on-reconstructi,ou2018augmented}, the system in \eqref{complex_system} is equivalent to the following system, which enforces the condition that $\mathbf{x}$ is real-valued:
\beq
\begin{pmatrix}
\rm{Re}(\mathbf{C})\\
\rm{Im}(\mathbf{C})\\
\mathbf{u}\\\mathbf{v}
\end{pmatrix}
\mathbf{x}=\begin{pmatrix}
\rm{Re}(\mathbf{d})\\
\rm{Im}({\mathbf{d}})\\
Q_{2N-1}\\0
\end{pmatrix}
\label{real_system}
\eeq
\item 
Solve \eqref{real_system} with an arbitrary-precision arithmetic system and then  apply the partial fraction decomposition to $[N-1/N](s)$ to find $\{  {r}_k, {p}_k\}_{k=1}^N$.
\end{enumerate}
Some numerical results from implementing this approach are shown in of Figures  \ref{fig:numericalE_N}  and \ref{compare}, where the typical behavior of $\epsilon(y)$ for various $N$ and $\alpha$ are demonstrated, with a clearly visible boundary layer on the left due to the branch cut of these functions at $z\le 0$. We report the numerical error $E_N$ with different $N$ and $\alpha$. Note  that the proposed method approximates $z^{-\alpha}$ very well and $E_N$ decreases as \(N\) increases. {In addition, the exponential decay nature is observed from Figure \ref{fig:numericalE_N} (b), which is consistent with  Theorem \ref{thm_N_conv}. }

%%%

%}

\section{Fast algorithm for the time domain Caputo derivative $D^\alpha$ \label{sec:algorithm}}
%%%%%%%%%%%%%
%\subsection{The unshifted scheme}
%\label{unshifted}
Based on \eqref{temp_Dalpha} and \eqref{approx_z_neg}, we derive the following 
\beq
\mathcal{L}[D^\alpha u](z)&=& z^\beta\, \hat{u}\cdot z,\quad \mbox{with} \quad \beta:=\alpha-1
\nonumber\\
&\approx& c_N (z \hat{u})+ \sum_{k=1}^N \frac{z \gamma_k}{z-\eta_k} \hat{u},\, c_N\in \field{R},\,\gamma_k>0,\, \beta_k<0 \nonumber\\
&=& c_N\left( z\hat{u}\right)+\left(\sum_{k=1}^N \gamma_k +\sum_{k=1}^N \frac{\gamma_k\eta_k}{z-\eta_k} \right) \hat{u}\nonumber\\
&=:& c_N\left( z \hat{u}\right)+\left(D_N +\sum_{k=1}^N \frac{\sigma_k}{z-\eta_k} \right) \hat{u}
\label{pole-res-no-shift}
\eeq

Recall that the coefficients in \eqref{pole-res-no-shift} are related to the pole-residue approximation of $Q(s)$ in \eqref{approx_Q} as follows
\beq
c_N=1-\sum_{k=1}^N \frac{ {r}_k}{ {p}_k};\quad \gamma_k=\frac{ r_k}{ {p}_k^2}; \quad D_N=\sum_{k=1}^N \gamma_k ; \quad \eta_k=1-\frac{1}{ {p}_k}
\label{pole-res-para}
\eeq

Applying the inverse Laplace transform on both sides of \eqref{pole-res-no-shift} leads to the following approximation of $D^\alpha u(t)$ for $t>0$
\beq
D^\alpha u(t)\approx c_N u'(t)+ D_N\, u(t)+  \sum_{k=1}^N  \sigma_k \int_0^t e^{\eta_k (t-\tau)} u(\tau)\, d\tau
\label{approx_time_domain}
\eeq

Values  of $c_N$ and $D_N$ in \eqref{approx_time_domain} for  various values of $\alpha$ are listed in Table \ref{table:cndn}. {Equation \eqref{approx_time_domain} reveals that the Caputo derivative can be considered as the sum of a small part of the standard derivative $u'$, a dissipation term, and a group of relaxation terms with exponentially decaying kernels.}

\begin{table}[h!]
  \begin{center}
\caption{\label{table:cndn}$c_N, D_N$ for  $D^\alpha$ with various values of $\alpha$, $Y=5$, $N=30$}
 \begin{tabular}{|c|c|c|} 
\hline
$\alpha$ & $c_N$& $D_N$\\
\hline
0.01 &5.103258770056953e-05& 1.051911296432398\\
\hline
0.1 & 8.336818977919296e-04& 1.645080161526623\\
\hline
0.3& 0.007290391807505& 4.210745812240900\\
\hline
0.5&0.034875160507518&  9.857397711438244\\
\hline
0.7&0.139977209093904& 19.662188927847950\\
\hline
0.9&0.522786577631653& 22.555809615609743\\
\hline
0.99999&0.999993522354222& 4.272681233022533e-04\\
\hline
\end{tabular}
  \end{center}
\end{table}
%\vspace{0.5in}

%%%%%%%%%%%%%%%%%%%%%%%%%%%%%%%%%%%%%%%%%%%%%%%%%%%%%%%%%%%%%%
\section{{Proposed method for solving the Caputo TFDs  \label{sec:proposed}
}}

Consider the initial value problem (IVP)
\begin{equation}
D^\alpha \tilde{u}(t) = G(t,\tilde{u}(t)), \qquad \tilde{u}(0)=u_0,
\label{ivp_org}
\end{equation}
in $(0,T)$, where $\alpha \in (0,1)$, 
\[
G : [0,T]\times \Pi \rightarrow \mathbb{R}^d,\quad T>0,\quad \Pi \subset \mathbb{R}^d \text{ open},
\]
and $D^\alpha$ is the Caputo $\alpha$-derivative defined in \eqref{Def:Caputo}.

Define  $u(t):=\tilde{u}(t)-u_0$ and $F(t,u(t)):=G(t,\tilde{u}(t)+u_0)$, then $D^\alpha u=D^\alpha \tilde{u}$ and $u(0)=0$. The IVP now reads:
\begin{equation}
D^\alpha {u}(t) = F(t,u(t)), \qquad {u}(0)=0,
\label{eq:ivp}
\end{equation}

We approximate the left-hand-side of \eqref{eq:ivp} by \eqref{approx_time_domain}, followed by introducing the auxiliary variables
\[
\phi_k(t) = \int_0^t e^{\eta_k (t-\tau)} u(\tau)\, d\tau.
\]
Applying the Leibniz rule to $\phi_k(t)$, we approximate \eqref{eq:ivp} with the following coupled system of equations for unknown functions \(u(t)\) and \(\phi_k(t)\):
\beq
c_N\frac{du}{dt} =
F(t,u(t)) - D_N u(t) - \displaystyle\sum_{k=1}^N \sigma_k \phi_k(t),
\qquad u(0)=0, \label{system_eq_1}\\
\frac{d\phi_k}{dt} = \eta_k \phi_k(t) + u(t), 
\qquad \phi_k(0)=0, \quad k=1,2,\dots,N.
\label{eq:auxphiu}
\eeq
where $c_N, D_N, \{\gamma_k,\eta_k\}_{k=1}^N$ are defined in \eqref{pole-res-para}.  In terms of the unknown vector
\[
\mathbf{y}(t) :=
\begin{bmatrix}
u(t) & \phi_1(t) & \phi_2(t) & \cdots & \phi_N(t)
\end{bmatrix}^t
\in \mathbb{R}^{N+1}.
\]

, the coupled system can be written in the following compact matrix form:
\beq
\frac{d\mathbf{y}}{dt}
=
\mathbf{A}\,\mathbf{y}(t)
+
\mathbf{F}(t,u(t))
\label{ODE_Matrix}
\eeq
where  $\mathbf{A} \in \mathbb{R}^{(N+1)\times(N+1)}$ and $\mathbf{F}\in \mathbb{R}^{N+1}$ are given by
\[
\mathbf{A} =
\begin{bmatrix}
-\dfrac{D_N}{c_N} & -\dfrac{\sigma_1}{c_N} & -\dfrac{\sigma_2}{c_N} & \cdots & -\dfrac{\sigma_N}{c_N} \\
1 & \eta_1 & 0 & \cdots & 0 \\
1 & 0 & \eta_2 & \cdots & 0 \\
\vdots & \vdots & \vdots & \ddots & \vdots \\
1 & 0 & 0 & \cdots & \eta_N
\end{bmatrix},\quad \mathbf{F}(t,u(t)) =
\begin{bmatrix}
\dfrac{F(t,u(t))}{c_N} \\ 0 \\0 \\ \vdots \\ 0
\end{bmatrix}
\]
with homogeneous initial condition $\mathbf{y}(0)=\mathbf{0}$. This formulation transforms the original fractional differential equation into an equivalent finite-dimensional system of ordinary differential equations, where the memory effects are captured through the auxiliary variables.

{Due to the almost diagonal structure of matrix $\mathbf{A}$, we propose to apply the Schur complement method for speeding up the numerical computation as follows in the context of the standard BDF2 discretization of \eqref{ODE_Matrix},  
 where the solution at $t=h$ is obtained by the standard trapezoidal rule. More precisely, we have 
\begin{equation}
\label{eq:BDF2}
\frac{3\mathbf{y}^n-4\mathbf{y}^{n-1}+\mathbf{y}^{n-2}}{2h}=A\mathbf{y}^n+F(t_n,u^n),
\end{equation}
which is equivalent to 
\begin{equation*}
\mathcal{A}\mathbf{y}^n=r^n,\ \mathcal{A}=\begin{bmatrix}
a & b^T\\
\ell & D
\end{bmatrix},\ \mathbf{y}^n=\begin{bmatrix}u^n\\  \phi^n\end{bmatrix}, \ r^n=\begin{bmatrix}r_u^n\\  r_\phi^n\end{bmatrix}.
\end{equation*}
In the above, we have denoted  $\lambda = \frac{3}{2h}$ and 
\begin{equation*}
 a=\lambda+\frac{D_N}{c_N},
 \ 
 b_k=\frac{\sigma_k}{c_N},
 \ 
 \ell=\begin{bmatrix}-1\\ \vdots\\ -1\end{bmatrix},
 \ 
 D=\begin{bmatrix}
\lambda-\eta_1 & 0 & \cdots & 0 \\
0 & \lambda-\eta_2 & \cdots & 0 \\
\vdots & \vdots & \ddots & \vdots \\
0 & 0 & \cdots &\lambda- \eta_N
\end{bmatrix}.
\end{equation*}
It is known that the corresponding Schur complement is  
\begin{equation*}
 S=a-b^TD^{-1}\ell=\lambda+\frac{D_N}{c_N}+
 \sum_{k=1}^N\frac{\sigma_k/c_N}{\lambda-\eta_k}.
\end{equation*}
Therefore, the solution for (\ref{eq:BDF2}) can be represented as 
\begin{equation}\label{eq:BDF2solution}
u^n=S^{-1}\left(r_u^n-b^T D^{-1}r^n_\phi\right),\quad
\phi^n =D^{-1}(r^n_\phi-\ell u^n).
\end{equation}
}

{In the case of linear problems, \eqref{ODE_Matrix} can also be solved exactly without discretization errors by applying  the block matrix exponential to the augmented system:
\[
\frac{d}{dt}
\begin{bmatrix}
\mathbf{y}(t)\\[2mm]
1
\end{bmatrix}
=
\begin{bmatrix}
\mathbf{A} & \mathbf{F}_n\\[1mm]
\mathbf{0} & 0
\end{bmatrix}
\begin{bmatrix}
\mathbf{y}(t)\\[2mm]
1
\end{bmatrix}.
\]
to obtain the following time-stepping scheme, 
\[
\begin{bmatrix}
\mathbf{y}_{n+1}\\[2mm]
1
\end{bmatrix}
=
\exp\!\left(
h
\begin{bmatrix}
\mathbf{A} & \mathbf{F}_n\\[1mm]
\mathbf{0} & 0
\end{bmatrix}
\right)
\begin{bmatrix}
\mathbf{y}_n\\[2mm]
1
\end{bmatrix},
\]
, which is referred to as the Block Exponential Integrator approach. The only numerical error arises from the computation of the matrix exponential in MATLAB. }  

After the solution is obtained, the solution to \eqref{ivp_org} is recovered as
\[
\tilde{u}(t)=u(t) + u_0.
\]
\section{{Numerical examples and comparison with other methods
}}
In this section, we demonstrate the performance of the numerical scheme proposed in Section \ref{sec:proposed} through four illustrative examples.
The first three problems are simple linear equations, namely one-dimensional fractional ordinary differential equations with constant coefficients, whose analytical solutions are known.
The first example exhibits growth behavior and the second one represents a decaying solution expressed in terms of the Mittag--Leffler function. In the third example, we consider the fractional harmonic oscillator as a test problem to evaluate the ability of the proposed method to capture the energy-preserving behavior of the system. {The fourth example is the fractional Van der Pol equation, which is nonlinear.} These examples enable us to investigate the accuracy and convergence properties of the proposed method.

{Moreover, we compare the results by the proposed method with those by two state-of-the-art numerical methods for solving TFDs: the $L_1$ method \cite{jin2023numerical}, which is a direct method, and the SOE method \cite{jiang2017fast}. }
{{
In the $L_1$ method,  the Caputo fractional derivative $D_t^\alpha y(t_n)$ is approximated  by \cite{jin2023numerical}
$$
D_t^\alpha y(t_n)
\approx \tau^{-\alpha} \left[ b_0 y(t_n) - b_{n-1} y(t_0) + \sum_{j=1}^{n-1} (b_j - b_{j-1}) y(t_{n-j}) \right],
$$
where $b_j$ are given by
$$
b_j=\frac{(j+1)^{1-\alpha}-j^{1-\alpha}}{\Gamma(2-\alpha)},\quad j=0,1,...,n-1.
$$
In the SOE method the approximation of   $D_t^\alpha y(t_n)$ is \cite{jiang2017fast}
\begin{equation*}
 D_t^\alpha y(t_n)
\approx
\frac{y^n-y^{n-1}}{\Delta t^\alpha \Gamma(2-\alpha)}
+
\frac{1}{\Gamma(1-\alpha)}
\left[
\frac{y^{n-1}}{\Delta t^\alpha}
-
\frac{y^0}{t_n^\alpha}
-
\alpha \sum_{i=1}^{N_{\exp}} \omega_i \mathcal{Y}_{\mathrm{h}}^i(t_n)
\right].
\end{equation*}
}}

\subsection{ Fractional Constant Coefficient ODE : Problem 1}

The first numerical test problem is 
 \beq
 {}^{C}D_t^{\alpha} y= 1 \,\,\,\,\,\,\,\,\,\,\,\,  y(0)=0, \quad 0<t<T
 \label{eq:exampleConstant}
 \eeq
 , whose exact solution is given by
\[ y=  \frac{t^{\alpha}}{\Gamma(\alpha +1)} .  \]
%%%%%figure 4%%%%%%%%%%%%%%%%%%%%%%%%%%%%%%%%%%%%
\begin{figure}[htbp]
\centering
\begin{subfigure}{0.48\linewidth}
    \centering
    \includegraphics[width=\linewidth]{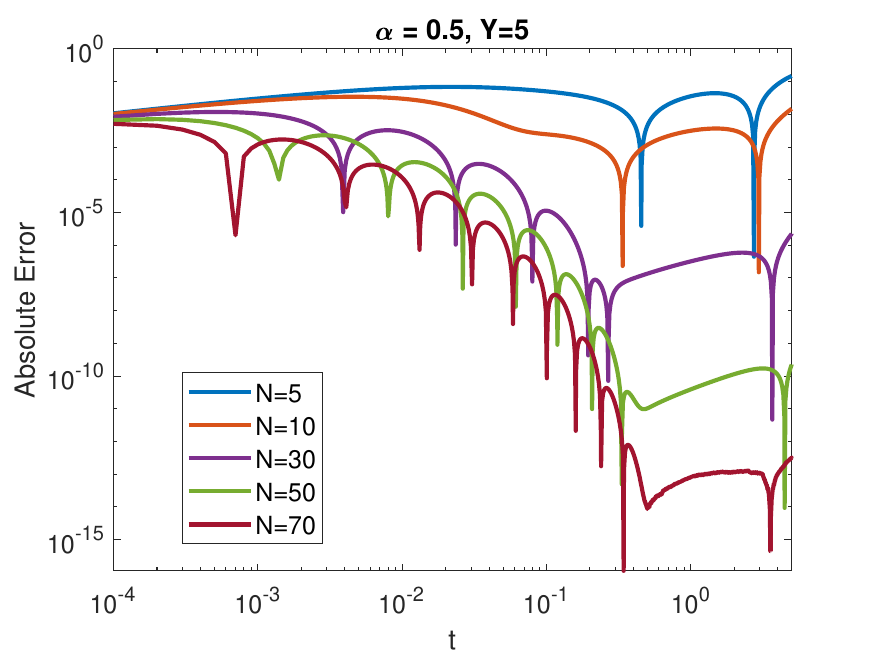}
    \caption{Local error on semi-log scale}
\end{subfigure}
\hfill
\begin{subfigure}{0.48\linewidth}
    \centering
    \includegraphics[width=\linewidth]{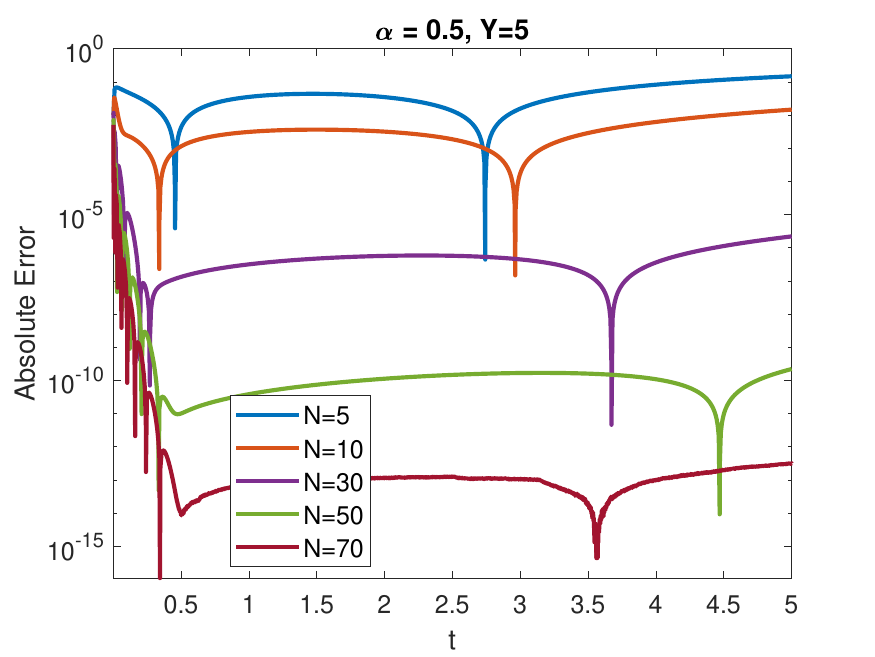}
    \caption{Local error on  log-log scale}
\end{subfigure}

\caption{Comparison of the absolute error for Problem~1. The left and right panels show the local error as a function of time and the local error on a logarithmic scale, respectively, for different numbers of quadrature points, $N=5,10,30,50,70$, with $\alpha=0.5$, $Y=5$, and $T=5$.}
\label{fig12_problem1}

\end{figure}

Figure \ref{fig12_problem1} compares the numerical solution obtained using the proposed block exponential integrator with the exact solution of Problem 1 for different $N$ in the pole-residue approximation  $N=5,10,30,50,70$, with $\alpha= 0.5$ and $Y=5$. As is shown in the figure, the pointwise error decreases systematically as $N$ increases, demonstrating the improved accuracy of the rational approximation. The left panel presents the local error over the entire time interval on a semilogarithmic scale, illustrating the overall error behavior. In contrast, the right panel shows the same error on a log--log scale, which provides a clearer view of the error evolution during the early stage of the simulation, where the effect of the weak singularity of the Caputo kernel is most pronounced.
%%%%%%%%%%%%%%%%%%%%%%%%

%%%%%%%%%Table 2%%%%%%%%%%%
 \begin{table}[http]
\centering
\caption{Problem 1. Absolute errors at selected time instants ($N=70$, $Y=5$, $\alpha=0.5$). The first and last ten values are reported, while the intermediate values are omitted for brevity.}
\label{tab:abs_error_short}
\begin{tabular}{cc}
\hline
Time & Absolute Error \\
\hline
0.00 & $5.1484\times10^{-3}$ \\
0.01 & $1.1106\times10^{-4}$ \\
0.02 & $3.5806\times10^{-5}$ \\
0.03 & $6.2832\times10^{-7}$ \\
0.04 & $4.6574\times10^{-6}$ \\
0.05 & $1.5975\times10^{-6}$ \\
0.06 & $1.4654\times10^{-7}$ \\
0.07 & $4.5063\times10^{-7}$ \\
0.08 & $2.8273\times10^{-7}$ \\
0.09 & $9.9972\times10^{-8}$ \\
\hline
$\vdots$ & $\vdots$ \\
\hline
4.91 & $3.0242\times10^{-13}$ \\
4.92 & $3.0242\times10^{-13}$ \\
4.93 & $2.9976\times10^{-13}$ \\
4.94 & $3.0420\times10^{-13}$ \\
4.95 & $3.0687\times10^{-13}$ \\
4.96 & $3.0598\times10^{-13}$ \\
4.97 & $3.0997\times10^{-13}$ \\
4.98 & $3.0997\times10^{-13}$ \\
4.99 & $3.1264\times10^{-13}$ \\
5.00 & $3.1930\times10^{-13}$ \\
\hline
\end{tabular}
\end{table}

Table \ref{tab:abs_error_short} presents the absolute (local) error at selected time instants during the early stage of the simulation. It is observed that the local error is relatively larger immediately after the initial time and then decreases rapidly as time progresses, {as is predicted by \eqref{error_estimate}}.

%VAR Y%%%%figure 5%%%%%%%%%%%%%%%%%%%%%%%%%%%%%%%%%%%%
\begin{figure}[htbp]
\centering
\begin{subfigure}{0.48\linewidth}
    \centering
    \includegraphics[width=\linewidth]{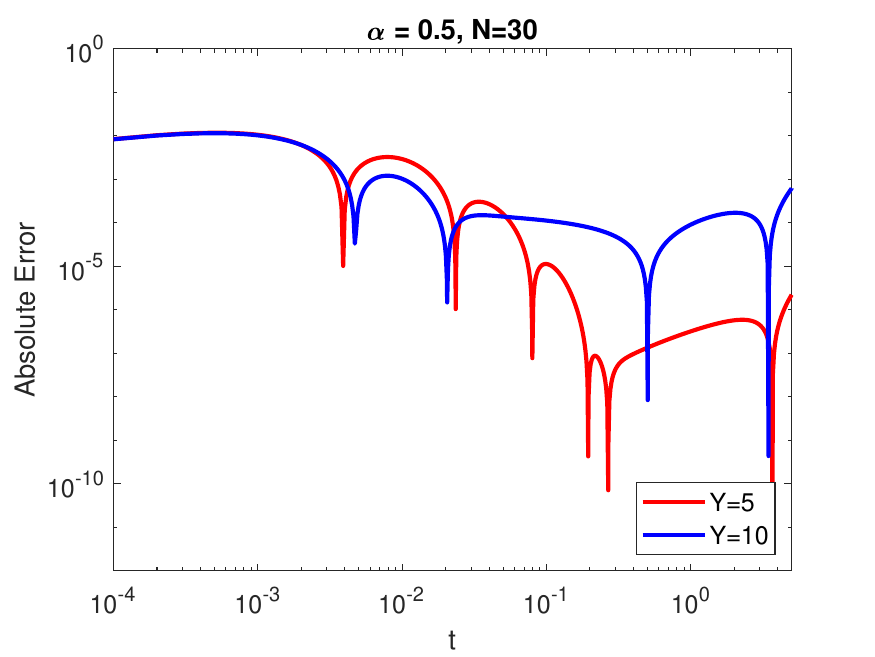}
    \caption{Local error on semi-log scale}
\end{subfigure}
\hfill
\begin{subfigure}{0.48\linewidth}
    \centering
    \includegraphics[width=\linewidth]{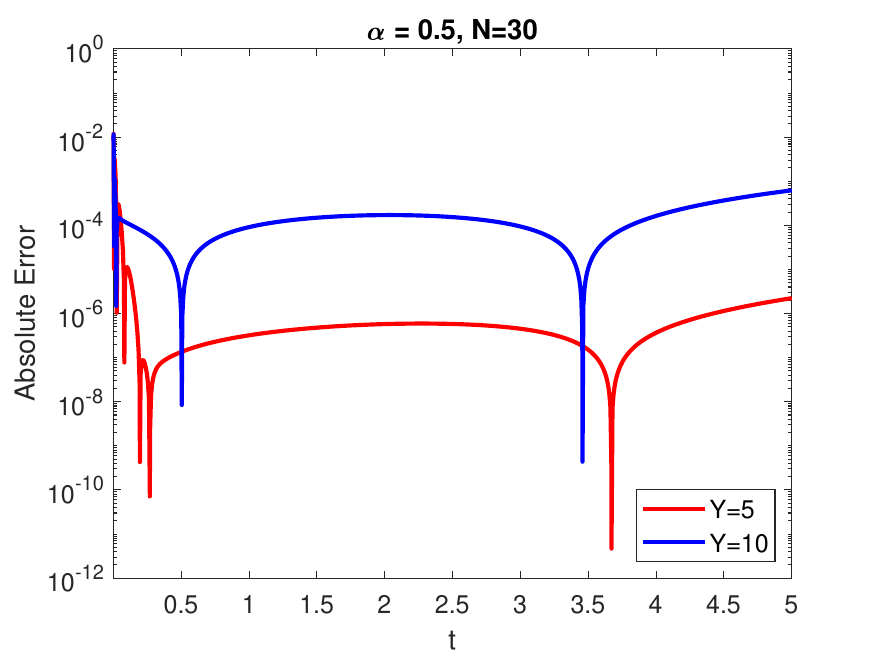}
    \caption{Local error on  log-log scale}
\end{subfigure}

\caption{Comparison of the absolute error for Problem~1 for different interpolation ranges \(Y\) with \(\alpha=0.5\), \(N=30\), and \(T=5\). The left and right panels show the local error as a function of time on semilogarithmic and log--log scales, respectively. The log--log plot highlights the influence of the interpolation range during the early stage of the simulation.}
\label{varY_problem1}

\end{figure}

%%%%EXPLANATION %%%%%%%%%%%%%%%%%%%%%%%%%%%%%%%%

Figure~\ref{varY_problem1} illustrates the effect of the interpolation range \(Y\) on the local error. For a fixed number of quadrature points \(N\), a smaller value of \(Y\) results in a smaller interpolation error \(\epsilon(y)\) in \eqref{error_estimate}, but a larger extrapolation error. As predicted by the analysis in Section~\ref{error-link}, the extrapolation error decays exponentially as time increases. Consequently, the error associated with the smaller value of \(Y\) eventually becomes lower than that corresponding to the larger value of \(Y\). This transition is observed around \(t\approx10^{-2}\), where the red curve crosses below the blue curve.

Although the choice of $Y$ has only a minor influence on the short-time interval, it becomes more significant over longer times. The results show that smaller point wise errors are obtained when a shorter interpolation range is used.
%%%ALL ALPHA%%%%%%%Figure 6%%%%%%%%%%%%%%%%%%%%%%%%%%%%%

\begin{figure}[htbp]
\centering
\begin{subfigure}{0.48\linewidth}
    \centering
    \includegraphics[width=\linewidth]{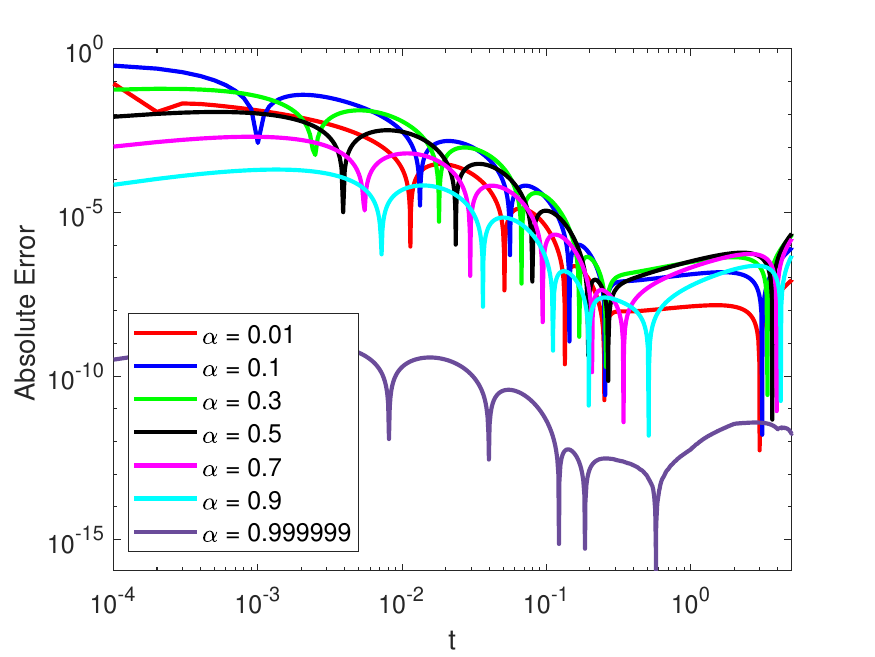}
    \caption{Local error on log-log scale}
\end{subfigure}
\hfill
\begin{subfigure}{0.48\linewidth}
    \centering
    \includegraphics[width=\linewidth]{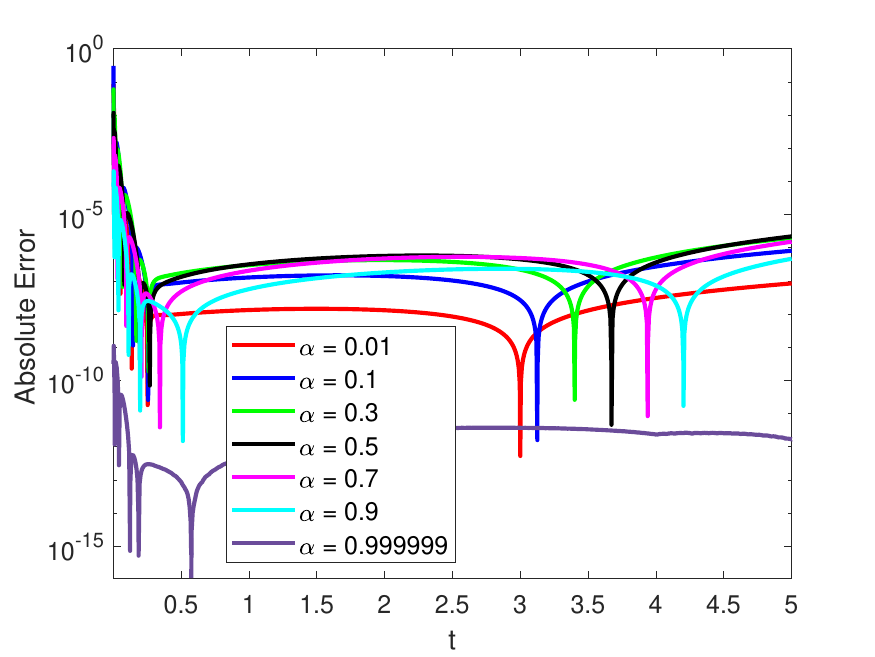}
    \caption{Local error on  semi-log scale}
\end{subfigure}

\caption{Comparison of the absolute error for Problem~1 for different values of \(\alpha\), with \(N=30\), \(Y=5\), and \(T=5\). The left and right panels show the local error as a function of time on log--log scales and semilogarithmic scales, respectively. The log--log plot provides a detailed view of the error behavior during the early stage of the simulation, while the semilogarithmic plot illustrates the overall error evolution over the entire time interval.}
\label{varalpha_problem1}

\end{figure}
%%%%%%explanation%%%%%%%%%%%%%%%%%

Figure~\ref{varalpha_problem1} presents a comparison of the absolute error for Problem~1 for different values of \(\alpha\), with \(N=30\), \(Y=5\), and \(T=5\). The left and right panels display the local error on log--log and semilogarithmic scales, respectively. The results show that the smallest error is obtained when \(\alpha\) is close to \(1\). Furthermore, the error exhibits a behavior near \(t=0\) similar to that observed over the entire time interval. Specifically, the error decreases rapidly during the initial stage of the simulation and subsequently stabilizes for larger times.

{\subsubsection{Comparison with the $L_1$ scheme and the SOE method}}
{To test the efficiency of the proposed method, we applied the BDF2 discretization and set $h = T/M$ with $M=50000,100000,...,250000$ to solve (\ref{eq:exampleConstant}) with $\alpha=0.5,\, Y=5$, $N=70$ and $T=5$ by all the three methods.} 

{Table \ref{Tabelconstant} consists of the the $L^2$ errors and orders of all the three methods.  Numerical results show that the convergence rate of the proposed methods is 1, highest among all the three methods. 
}

\begin{figure}[htbp]
\centering
    \begin{subfigure}{0.48\linewidth}
        \centering
        \includegraphics[width=\linewidth]{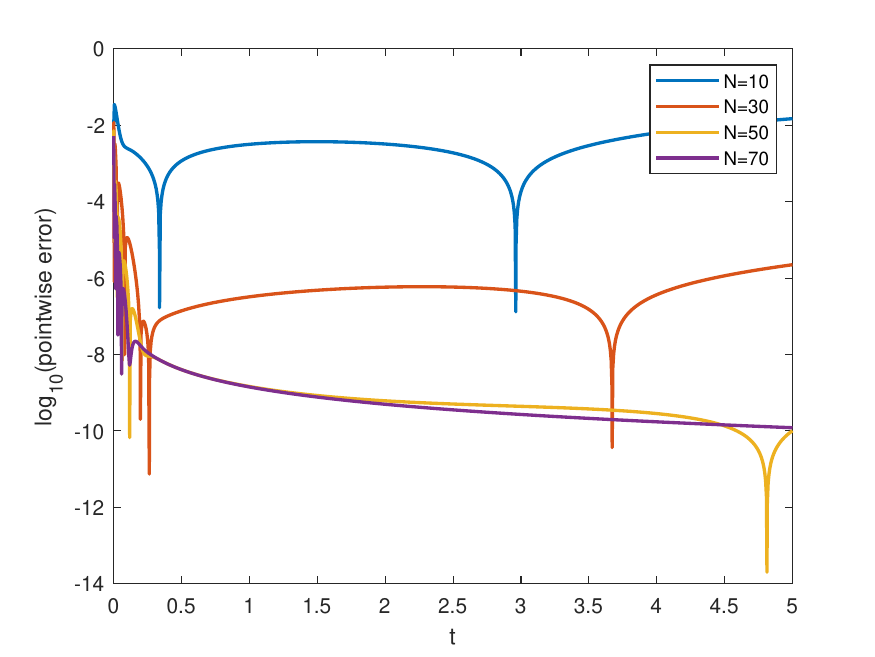}
        \caption{}
    \end{subfigure}\hfill
    \begin{subfigure}{0.48\linewidth}
        \centering
        \includegraphics[width=\linewidth]{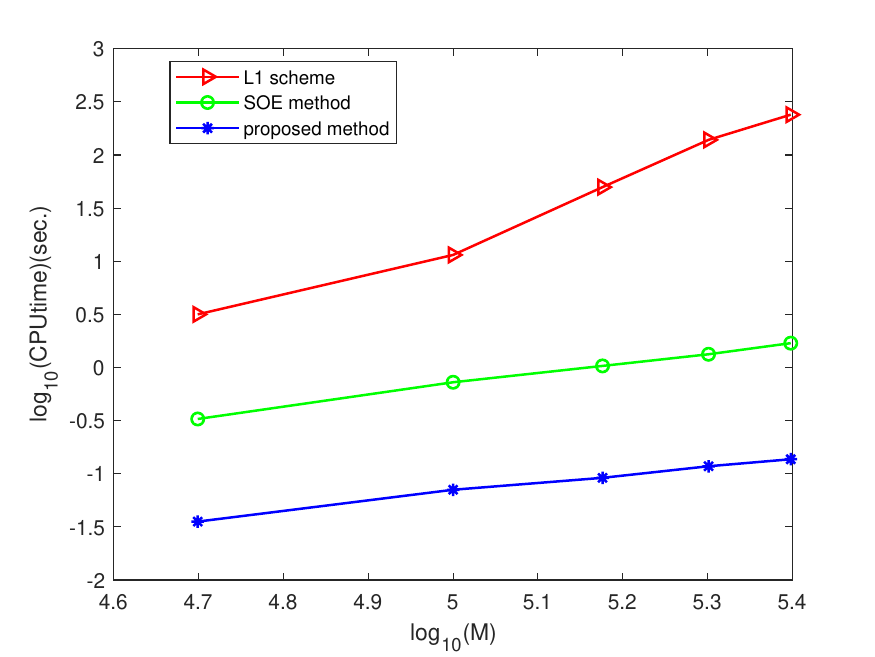}
        \caption{}
    \end{subfigure}
    \caption{Numerical results for Problem 1. (a) the point-wise errors for fixed $\alpha=0.5$ and $M=5\times 10^4$ with various $N=10,30,50,70$, (b) the CPU time with $Y=5,\,\alpha=0.5$ and $N=70$. \label{fig:pointerrorsCPU}}
\end{figure}

{In Figure \ref{fig:pointerrorsCPU} (a), we report the point-wise error for different $N$. Note  that the error decreases with increasing $N$ but gets saturated near $N=50$ because the discretization error corresponding to the fixed $M$ overshadows the small error resulting from a larger $N$.
}  

{ Figure \ref{fig:pointerrorsCPU} (b) shows the linear complexity of the proposed method, which also has the lowest CPU time among the three methods.
}
\begin{table}
\centering
{
\begin{tabular}{p{10pt} p{5pt} p{50pt} p{20pt}  p{5pt} p{50pt} p{20pt} p{5pt} p{50pt} p{20pt} }
\hline
   \multirow{2}{*}{M}& & \multicolumn{2}{c}{$L_1\, method$}&  & \multicolumn{2}{c}{$SOE\, method$}& & \multicolumn{2}{c}{$Proposed\, method$} \\
\cline{3-4} \cline{6-7} \cline{9-10}
  & & Error &Order     & &    Error &Order &    & Error &Order \\
\hline
 $100$  &  & 3.107E-2  & -        & &3.107E-2  & -       &    & 1.531E-2   & - \\
 $200$  &  & 1.660E-2  &  0.904   & &1.660E-2  &  0.904  &    & 7.657E-3   & 1.000\\
 $400$  &  & 8.805E-3  &  0.915   & &8.805E-3  &  0.915  &    & 3.828E-3   & 1.000 \\
 $800$  &  & 4.641E-3  &  0.924   & &4.641E-3  &  0.924  &    & 1.901E-3   & 1.004 \\
 $1600$ &  & 2.433E-3  &  0.932   & &2.433E-3  &  0.932  &    & 9.343E-4   & 1.030 \\
\hline
\end{tabular}
}
\caption{The $L^2$ errors and rates for Problem 1 with $\alpha = 0.5,\,Y=5,\,N=70$.}\label{Tabelconstant}
\end{table}

\subsection{Fractional Constant Coefficient ODE: Problem 2 }
The second  numerical test problem is

\begin{equation}
D^{\alpha} u = -u, \qquad u(0) = 1, \qquad t \in (0,T).
\label{example2}
\end{equation}
The  exact solution of \eqref{example2} is given in terms of the Mittag-Leffler function 
\begin{equation}
u(t) = E_{\alpha}(-t^{\alpha}):=\sum_{k=0}^{\infty} \frac{t^{k}}{\Gamma(\alpha k + 1)}.
\end{equation}
We first examine the effect of the number of quadrature points on the absolute error. Table~\ref{tab:N30_N100} presents the absolute errors at selected time instants for different numbers of quadrature points, \(N=30,40,50,60,70,\) and \(100\), with \(\alpha=0.5\), \(Y=5\), and \(T=5\). As can be seen from the table, the local error decreases significantly as the number of quadrature points increases, demonstrating the improved accuracy of the proposed approximation. In particular, the errors are reduced from approximately \(10^{-7}\) for \(N=30\) to nearly machine precision for \(N=100\) over most of the computational interval. Moreover, the CPU times remain relatively small for all values of \(N\), indicating that increasing the number of quadrature points has only a modest impact on the overall computational cost while providing a substantial gain in accuracy.

%%%%%%%%%%%%%%%%%%%%%%%%%%%%%%%%%
\begin{table}[http]
\centering
\caption{Problem 2: Comparison of absolute errors at selected time instants for different numbers of quadrature points (\(\alpha=0.5\), \(Y=5\), and \(T=5\)).}
\label{tab:N30_N100}
\begin{adjustbox}{width=\textwidth}
\begin{tabular}{cccccccc}
\hline
Time & $N=30$ & $N=40$ & $N=50$ & $N=60$ & $N=70$ & $N=100$ \\
\hline
0.0 & $0.0000$ & $0.0000$ & $0.0000$ & $0.0000$ & $0.0000$ & $0.0000$ \\

0.5 & $8.6604\times10^{-8}$ & $8.2267\times10^{-10}$ &
$8.1848\times10^{-12}$ & $1.2468\times10^{-13}$ &
$9.1038\times10^{-15}$ & $3.3307\times10^{-16}$ \\

1.0 & $1.1454\times10^{-7}$ & $1.3167\times10^{-9}$ &
$1.5653\times10^{-11}$ & $1.7975\times10^{-13}$ &
$2.7756\times10^{-15}$ & $2.1094\times10^{-15}$ \\

1.5 & $1.2348\times10^{-7}$ & $1.6763\times10^{-9}$ &
$2.2791\times10^{-11}$ & $3.0453\times10^{-13}$ &
$1.0547\times10^{-14}$ & $2.6645\times10^{-15}$ \\

2.0 & $1.0927\times10^{-7}$ & $1.8054\times10^{-9}$ &
$2.8225\times10^{-11}$ & $4.2094\times10^{-13}$ &
$2.4591\times10^{-14}$ & $5.3846\times10^{-15}$ \\

2.5 & $6.7790\times10^{-8}$ & $1.5939\times10^{-9}$ &
$3.0082\times10^{-11}$ & $5.2636\times10^{-13}$ &
$2.9643\times10^{-14}$ & $8.3267\times10^{-15}$ \\

3.0 & $5.0187\times10^{-9}$ & $9.2008\times10^{-10}$ &
$2.6033\times10^{-11}$ & $5.6327\times10^{-13}$ &
$4.0801\times10^{-14}$ & $1.7708\times10^{-14}$ \\

3.5 & $1.1318\times10^{-7}$ & $3.4918\times10^{-10}$ &
$1.3018\times10^{-11}$ & $2.3126\times10^{-13}$ &
$2.8388\times10^{-13}$ & $2.7189\times10^{-13}$ \\

4.0 & $2.6061\times10^{-7}$ & $2.3647\times10^{-9}$ &
$1.9583\times10^{-11}$ & $7.8155\times10^{-12}$ &
$8.0653\times10^{-12}$ & $8.0602\times10^{-12}$ \\

4.5 & $4.5125\times10^{-7}$ & $5.4258\times10^{-9}$ &
$2.2141\times10^{-10}$ & $1.6946\times10^{-10}$ &
$1.6919\times10^{-10}$ & $1.6919\times10^{-10}$ \\

5.0 & $6.9084\times10^{-7}$ & $1.1790\times10^{-8}$ &
$2.6907\times10^{-9}$ & $2.5786\times10^{-9}$ &
$2.5774\times10^{-9}$ & $2.5775\times10^{-9}$ \\

\hline
CPU Time (s) &
$1.023\times10^{-1}$ &
$1.3299\times10^{-1}$ &
$1.3852\times10^{-1}$ &
$1.8777\times10^{-1}$ &
{$2.0161\times10^{-1}$} &
$6.7055\times10^{-1}$ \\
\hline
\end{tabular}
\end{adjustbox}
\end{table}

\begin{table}[http]
\centering
\caption{Problem 2: Comparison of absolute errors at selected time instants for different values of $\alpha$ ($N=30$, $Y=5$, $T=5$). CPU Time = 1.3179 sec. for all computations. }
\label{tab:mittag_alpha_comparison}
\begin{adjustbox}{width=\textwidth}
\begin{tabular}{cccccccc}
\hline
$t$ &
$\alpha=0.01$ &
$\alpha=0.1$ &
$\alpha=0.3$ &
$\alpha=0.5$ &
$\alpha=0.7$ &
$\alpha=0.9$ &
$\alpha=0.99999$ \\
\hline
0.0 &
$0.0000$ &
$0.0000$ &
$0.0000$ &
$0.0000$ &
$0.0000$ &
$0.0000$ &
$0.0000$ \\

0.5 &
$2.5801\times10^{-9}$ &
$2.5254\times10^{-8}$ &
$6.7005\times10^{-8}$ &
$8.6604\times10^{-8}$ &
$7.4215\times10^{-8}$ &
$2.8323\times10^{-8}$ &
$2.1694\times10^{-13}$ \\

1.0 &
$3.3675\times10^{-9}$ &
$3.2990\times10^{-8}$ &
$8.7793\times10^{-8}$ &
$1.1454\times10^{-7}$ &
$1.0179\times10^{-7}$ &
$4.4069\times10^{-8}$ &
$2.8477\times10^{-13}$ \\

1.5 &
$3.6462\times10^{-9}$ &
$3.5725\times10^{-8}$ &
$9.4934\times10^{-8}$ &
$1.2348\times10^{-7}$ &
$1.0973\times10^{-7}$ &
$4.8472\times10^{-8}$ &
$2.1377\times10^{-13}$ \\

2.0 &
$3.2653\times10^{-9}$ &
$3.1989\times10^{-8}$ &
$8.4699\times10^{-8}$ &
$1.0927\times10^{-7}$ &
$9.5963\times10^{-8}$ &
$4.2020\times10^{-8}$ &
$4.6296\times10^{-14}$ \\

2.5 &
$2.0740\times10^{-9}$ &
$2.0320\times10^{-8}$ &
$5.3426\times10^{-8}$ &
$6.7790\times10^{-8}$ &
$5.8017\times10^{-8}$ &
$2.4649\times10^{-8}$ &
$2.2300\times10^{-13}$ \\

3.0 &
$7.9263\times10^{-11}$ &
$6.6772\times10^{-10}$ &
$2.5447\times10^{-9}$ &
$5.0187\times10^{-9}$ &
$6.6338\times10^{-9}$ &
$3.8960\times10^{-9}$ &
$5.8996\times10^{-13}$ \\

3.5 &
$3.3481\times10^{-9}$ &
$3.0368\times10^{-8}$ &
$8.6889\times10^{-8}$ &
$1.1318\times10^{-7}$ &
$1.0047\times10^{-7}$ &
$4.3949\times10^{-8}$ &
$1.0124\times10^{-12}$ \\

4.0 &
$7.8715\times10^{-9}$ &
$3.7680\times10^{-8}$ &
$2.0326\times10^{-7}$ &
$2.6061\times10^{-7}$ &
$2.2586\times10^{-7}$ &
$9.5850\times10^{-8}$ &
$1.7228\times10^{-13}$ \\

4.5 &
$1.3513\times10^{-8}$ &
$3.3540\times10^{-7}$ &
$3.5545\times10^{-7}$ &
$4.5125\times10^{-7}$ &
$3.8506\times10^{-7}$ &
$1.5995\times10^{-7}$ &
$3.1970\times10^{-11}$ \\

5.0 &
$1.6447\times10^{-8}$ &
$4.1178\times10^{-6}$ &
$5.4959\times10^{-7}$ &
$6.9084\times10^{-7}$ &
$5.8097\times10^{-7}$ &
$2.3698\times10^{-7}$ &
$5.7703\times10^{-10}$ \\
\hline
\end{tabular}
\end{adjustbox}

\end{table}
%%%%%%%%%%%%%%%%%%%%%%%%%%%%%%%%%%%%%%%%%%%%%%%

Finally, we investigate the effect of the fractional order \(\alpha\) on the local error. For this purpose, Table~\ref{tab:mittag_alpha_comparison} presents the absolute errors at selected time instants for different values of \(\alpha\), while keeping \(N=30\), \(Y=5\), and \(T=5\) fixed. The results show that the local error decreases as the fractional order \(\alpha\) approaches \(1\). In particular, the smallest absolute errors are obtained for values of \(\alpha\) close to \(1\), indicating that the proposed method becomes increasingly accurate as the fractional-order model approaches the corresponding classical first-order differential equation.
%%%%%%%%%%Comparison%%%%%%%%%%%%%%%%%%%%%%%%%%%%%%%%%%%%%%%%%%%%%%%

{\subsubsection{Comparison with the $L_1$ method and the SOE method}}
{
To test the efficiency of the proposed method, we applied the BDF2 discretization and set $h=T/M$ with $T=10$ and $M=40000,\,80000,\,...,\,200000$.}
{Table \ref{TabelProblem2} 
shows the  $L^2$ errors and rates of the three method. The parameters of the proposed method are $\alpha=0.5,\,Y=5,\,N=70$. Similar to Problem 1, the proposed method has the highest convergence order than the other two methods.}

\begin{figure}[htbp]
\centering
    \begin{subfigure}{0.48\linewidth}
        \centering
        \includegraphics[width=\linewidth]{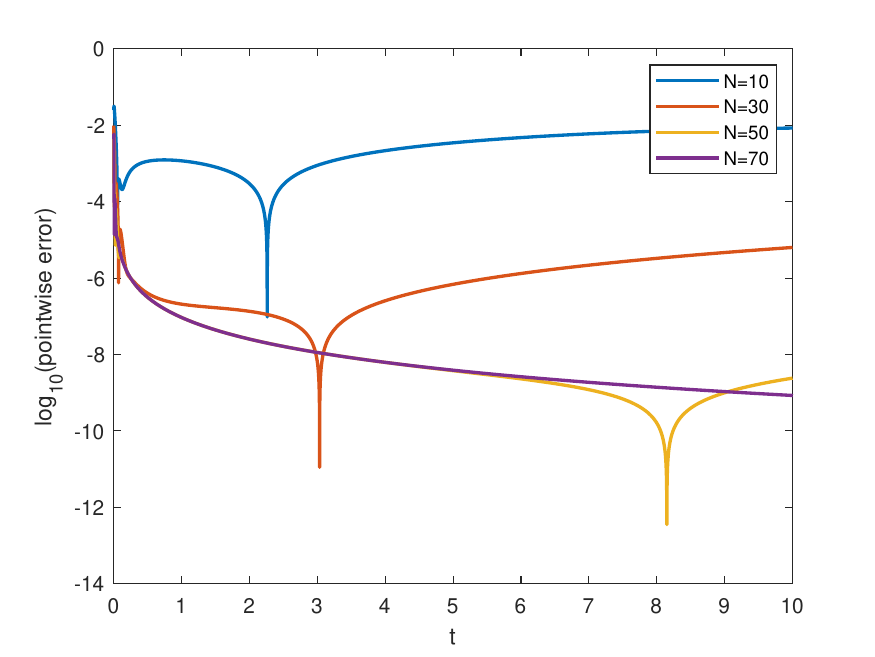}
        \caption{}
    \end{subfigure}\hfill
    \begin{subfigure}{0.48\linewidth}
        \centering
        \includegraphics[width=\linewidth]{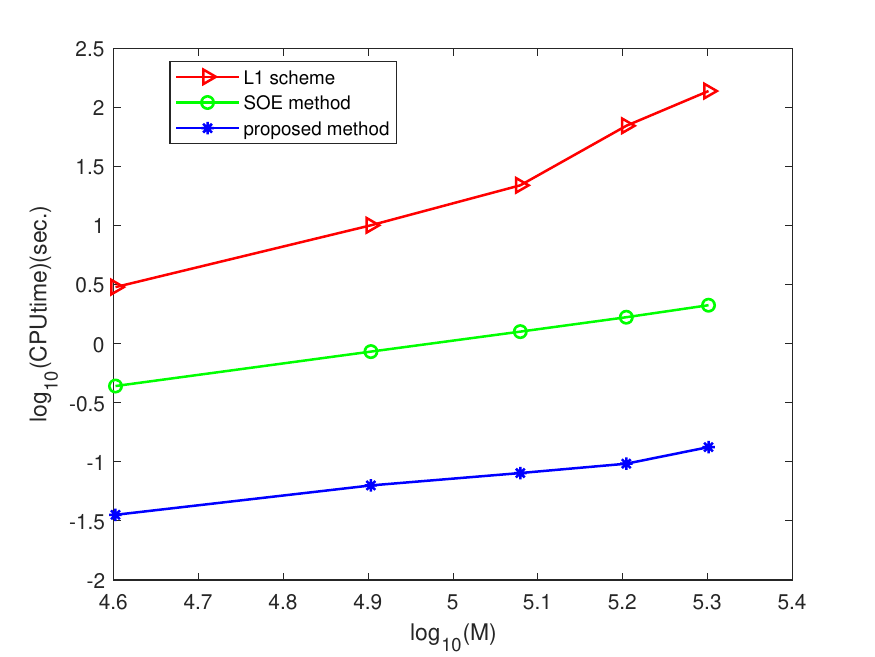}
        \caption{}
    \end{subfigure}
    \caption{Numerical results for Problem 2. (a) the point wise errors for fixed $\alpha=0.5$ and $M=10^4$ with various $N=10,30,50,70$, (b) the CPU time with $\alpha=0.5$ and $N=100$. \label{fig:pointerrorsCPUP2}}
\end{figure}

{In Figure \ref{fig:pointerrorsCPUP2} (a), we report the point-wise error for different $N$. Note  that the error decreases with increasing $N$.} 

{Figure \ref{fig:pointerrorsCPUP2} (b) shows the computational cost of the $L_1$, the SOE and the proposed methods. Here we have set $N=100$. The solution of SOE is obtained with $N_{exp}=74$. Compared with the $L_1$ method, linear complexity and a significant reduction in CPU time of the proposed method are observed. In addition, the proposed method is also less expensive than the SOE method due to the Schur complement.}

\begin{table}
\centering
{
\begin{tabular}{p{10pt} p{5pt} p{40pt} p{20pt}  p{5pt} p{40pt} p{20pt} p{5pt} p{40pt} p{20pt} }
\hline
   \multirow{2}{*}{M}& & \multicolumn{2}{c}{$L_1\, method$}&  & \multicolumn{2}{c}{$SOE\, method$}& & \multicolumn{2}{c}{$Proposed\, method$} \\
\cline{3-4} \cline{6-7} \cline{9-10}
  & & Error &Order     & &    Error &Order &    & Error &Order \\
\hline
 $1000$  &  & 4.20E-3    & -        & &4.20E-3  & -       &    & 3.05E-3   & - \\
 $2000$  &  & 2.33E-3    & 0.853    & &2.33E-3  & 0.853   &    & 1.52E-3   & 1.01\\
 $4000$  &  & 1.27E-3    & 0.870    & &1.27E-3  &  0.870  &    & 7.37E-4   & 1.04 \\
 $8000$  &  & 6.90E-4    & 0.884    & &6.90E-4  &  0.884  &    & 3.37E-4   & 1.13 \\
 $16000$ &  & 3.71E-4    & 0.896    & &3.71E-4  &  0.896  &    & 1.53E-4   & 1.14   \\
\hline
\end{tabular}
}
\caption{The $L^2$ errors and rates for Problem 2 with fixed $\alpha = 0.5$, $Y=5$ and $N=70$ and various $M$.}\label{TabelProblem2}
\end{table}

%
%%%%%%%%PROBLEM3%%%%%%%%%%%%%%%%%%%%%%%%%%%%%%%%%%%%%%%%%%%%%%%%%%
\subsection{Fractional Harmonic Oscillation: Problem 3}

In this example, we consider the IVP of fractional harmonic oscillation differential equation 

\begin{equation}
(D^{\alpha})^2 x(t) + x(t)=0,\quad x(0)=0 \mbox{ and } D^{\alpha} x(0)=1.
\label{harmonic}
\end{equation}
For $\alpha=1$, Equation~\eqref{harmonic} reduces to the classical harmonic oscillator equation. In the numerical experiments, we consider the fractional regime $0<\alpha<1$. According to the classical theory, the system is stable and energy-conservative when $\alpha=1$. Therefore, the proposed numerical method is expected to recover this property as $\alpha \to 1$. For $\alpha<1$, however, the energy is not preserved and the solution exhibits a spiral-type behaviour.

To apply the proposed methods, we rewrite  Equation~\eqref{harmonic} as  a system by substituting $y= D^{\alpha} x$:
\begin{equation}
D^{\alpha}x(t)=y(t),\qquad
D^{\alpha}y(t)=-x(t),
\qquad x(0)=0,\;y(0)=1.
\label{eq:fho}
\end{equation}

For the classical case $\alpha=1$, the exact solution is
\[
x(t)=\sin t,\qquad y(t)=\cos t,
\]
and the energy/invariant is preserved:
\begin{equation}
x(t)^2+y(t)^2 \equiv 1.
\label{eq:invariant}
\end{equation}

To test whether the proposed method conserve the energy as $\alpha\rightarrow 1$, we fix  $\alpha=0.999999$.

For $\alpha=0.9 <1$, the identity \eqref{eq:invariant} is generally not preserved; it is used here as a reference diagnostic when $\alpha$ is close to $1$.

In the numerical simulation, we fix $Y=5$, $N=30$  and $T=100$.
%%%%%figure 10%%%%%%%%%%%%%%%%%%%%%%%%%%%%%%%%%%%%
\begin{figure}[htbp]
\centering
\begin{subfigure}{0.48\linewidth}
    \centering
    \includegraphics[width=\linewidth]{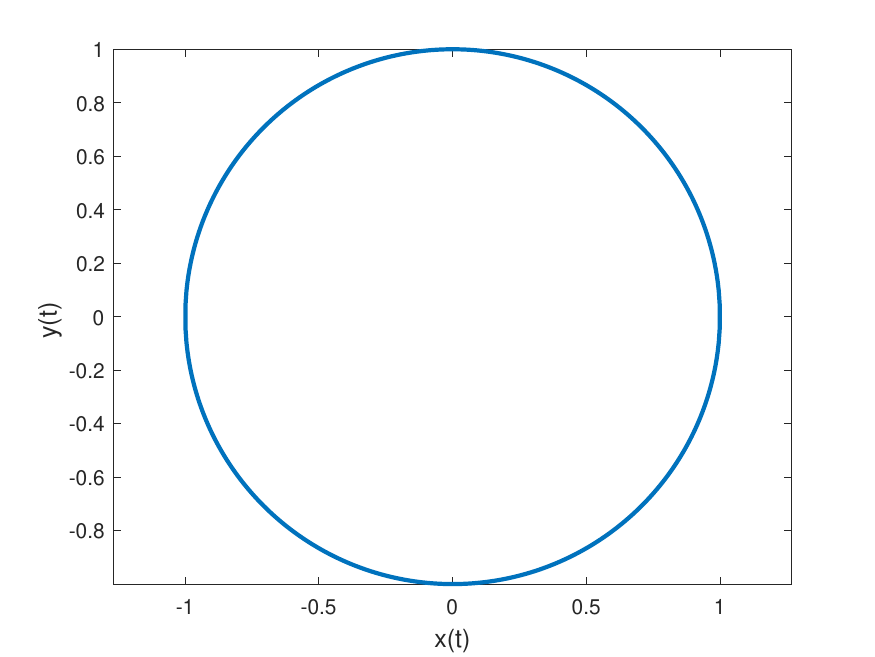}
    \caption{$\alpha=0.999999$.}
\end{subfigure}
\hfill
\begin{subfigure}{0.48\linewidth}
    \centering
    \includegraphics[width=\linewidth]{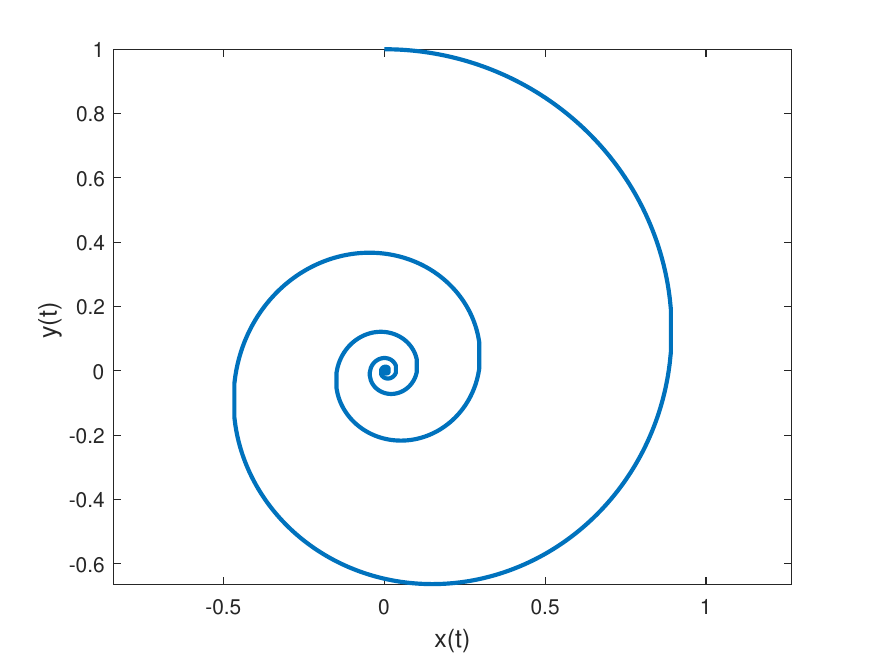}
    \caption{$\alpha=0.9$.}
\end{subfigure}

\caption{Problem 3. Numerical phase solutions of the fractional harmonic oscillation problem for $Y=5$, $N=30$, and $T=100$ for two different values of $\alpha=0.999999$ and $\alpha=0.9$}
\label{fig10_harmonic}

\end{figure}

Figure~\ref{fig10_harmonic} shows that the proposed method preserves the energy as $\alpha \to 1$. In this case, the fractional oscillation equation exhibits a closed orbit. However, when $\alpha=0.9<1$, the solution becomes spiral, indicating that the energy is no longer preserved.

\begin{figure}[htbp]
    \centering
        \includegraphics[width=0.7\linewidth]{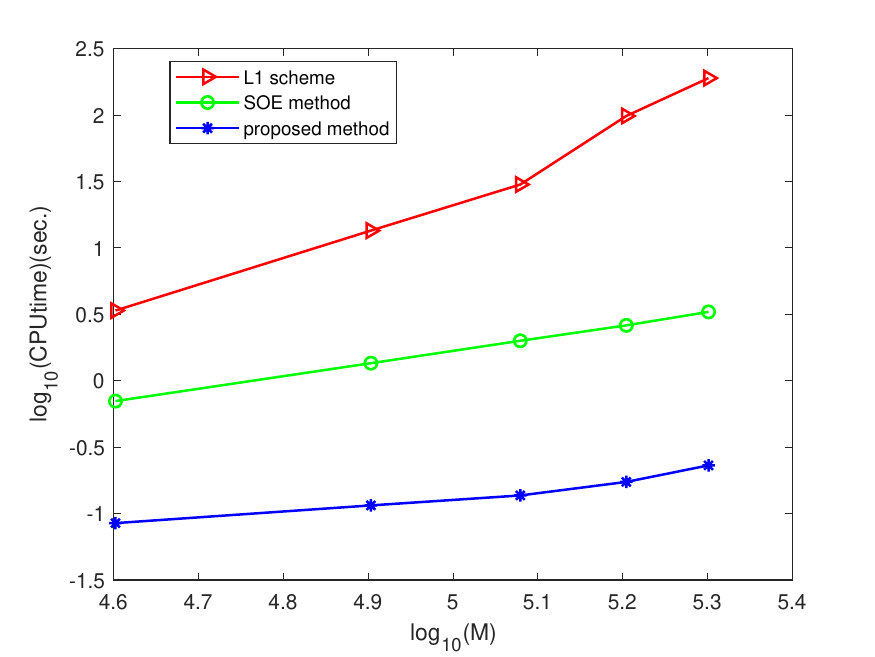}
            \caption{The CPU time for Problem 3 with  $\alpha=0.999999$, $Y=5$, $N=30$, $T=10$ and different $M$.}
        \label{fig:CPU_P3}
\end{figure}

\begin{figure}[htbp]
\centering
    \begin{subfigure}{0.48\linewidth}
        \centering
        \includegraphics[width=\linewidth]{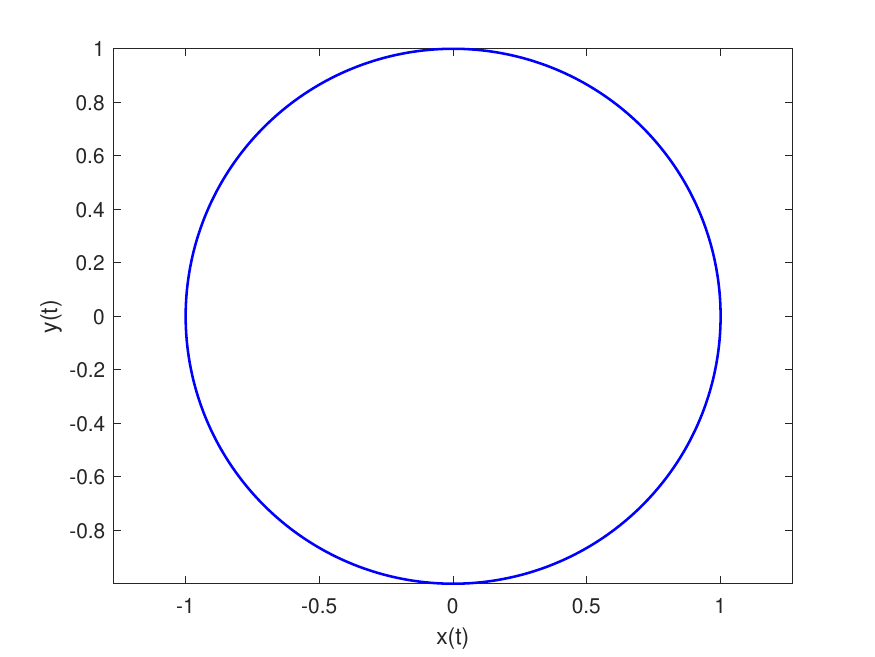}
        \caption{}
    \end{subfigure}\hfill
    \begin{subfigure}{0.48\linewidth}
        \centering
        \includegraphics[width=\linewidth]{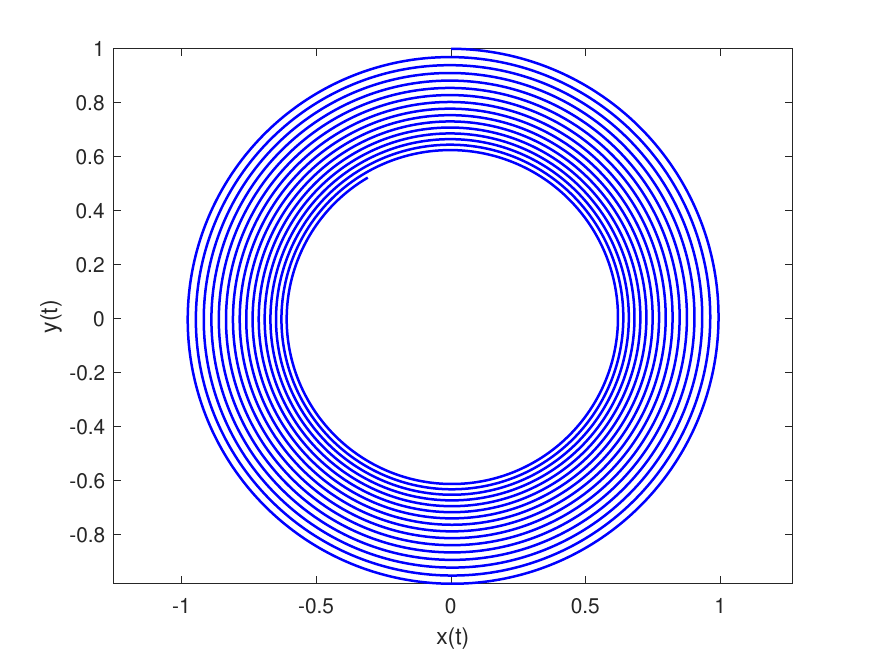}
        \caption{}
    \end{subfigure}
    \caption{{Numerical phase solutions of Problem 3 with $Y=5,\,N=30,\,\alpha=0.999999,\,T=100$ and $h=0.01$. (a) the numerical solution by the proposed method, (b) the numerical solution by the SOE method.
    \label{fig:PhaseP3}}}
\end{figure}

{\subsubsection{Comparison with the $L_1$ method and the SOE method} }
{
Setting $v=y-1$ and applying (\ref{system_eq_1}) - (\ref{eq:auxphiu}), we derive from (\ref{harmonic})
\begin{align}
 c_N \frac{dx}{dt} &=-D_Nx+v+1-\sum_{k=1}^{N}\sigma_k\phi_{x,k},
 \label{eq:ode_x}\\
 c_N \frac{dv}{dt} &=-x-D_Nv-\sum_{k=1}^{N}\sigma_k\phi_{v,k},
 \label{eq:ode_v}\\
 \frac{d \phi_{x,k}}{dt}&=\eta_k\phi_{x,k}+x,
 \label{eq:ode_phix}\\
 \frac{d \phi_{v,k}}{dt}&=\eta_k\phi_{v,k}+v,\ k=1,...,N.
 \label{eq:ode_phiv}
\end{align}
The BDF2 method for the above system reads  
\begin{equation} \label{eq:blockbdf_system}
\begin{bmatrix}
 \mathbf{A}_0 & \mathbf{B}\\
 \mathbf{E}   & \mathbf{C}
\end{bmatrix}
\begin{bmatrix}
 \mathbf{U}\\ \mathbf{W}
\end{bmatrix}
=
\begin{bmatrix}
 \mathbf{R_U}\\ \mathbf{R_W}
\end{bmatrix}, \quad
\mathbf{U}=\begin{bmatrix}x^n\\ v^n\end{bmatrix},\quad
\mathbf{W}=\begin{bmatrix}\phi_x^n\\ \phi_v^n\end{bmatrix},
\end{equation}
where we have denoted 
\begin{equation*}
 \mathbf{A}_0=
 \begin{bmatrix}
  a&-q\\ q&a
 \end{bmatrix},
 \qquad
 \mathbf{B}=
 \begin{bmatrix}
  \bm b^T&\bm 0^T\\
  \bm 0^T&\bm b^T
 \end{bmatrix},
\label{eq:A0_B}
\end{equation*}
\begin{equation*}
 \mathbf{E}=
 \begin{bmatrix}
  \ell&\bm 0\\
  \bm 0&\ell
 \end{bmatrix},
 \qquad
 \mathbf{C}=
 \begin{bmatrix}
  \bm D&\bm 0\\
  \bm 0&\bm D
 \end{bmatrix},
\label{eq:E_C}
\end{equation*}
with $q = \frac{1}{c_N}$ and $a$ given by the previous example. It is easy to derive the Schur complement is
$$
\mathbf{S} =\mathbf{A}_0-\mathbf{B}\mathbf{C}^{-1}\mathbf{E},
$$ 
and thus the solution is 
$$
\mathbf{U}= \mathbf{S}^{-1}\left( \mathbf{R_U}-\mathbf{B}\mathbf{C}^{-1}\mathbf{R_W} \right).
$$
}
{{
In Table \ref{TabelProblem3}, we demonstrates the $L^2$ errors and convergence rates, where we have set $Y=5$, $N=30$, $T=10$ and $h=T/M$ with different $M$. We can see that the proposed method achieves second-order convergence, whereas the $L_1$ and SOE methods exhibit first-order convergence. Figure \ref{fig:CPU_P3} displays the computational cost, which confirms that the proposed method achieves linear complexity and a significant reduction in CPU time.
In Figure \ref{fig:PhaseP3}, we illustrates the phase solution by the SOE and proposed methods, where we have set $T=100$ and $h=0.01$. Numerical results indicate that, as $\alpha\to 1$, the proposed method preserves the energy and is stable for long-time simulation, while the solution by the SOE method is over-dissipative in this case.
\begin{table}
\centering
{
\begin{tabular}{p{10pt} p{5pt} p{40pt} p{20pt}  p{5pt} p{40pt} p{20pt} p{5pt} p{40pt} p{20pt} }
\hline
   \multirow{2}{*}{M}& & \multicolumn{2}{c}{$L_1\, method$}&  & \multicolumn{2}{c}{$SOE\, method$}& & \multicolumn{2}{c}{$Proposed\, method$} \\
\cline{3-4} \cline{6-7} \cline{9-10}
  & & Error &Order     & &    Error &Order &    & Error &Order \\
\hline
 $1000$  &  & 8.96E-2    & -        & &8.96E-2  & -       &    & 6.08E-4   & - \\
 $2000$  &  & 4.52E-2    & 0.987    & &4.52E-2  &  0.987  &    & 1.52E-4   & 2.000\\
 $4000$  &  & 2.27E-2    & 0.993    & &2.27E-2  &  0.993  &    & 3.80E-5   & 2.000 \\
 $8000$  &  & 1.14E-2    & 0.994    & &1.14E-2  &  0.994  &    & 9.51E-6   & 2.000 \\
 $16000$ &  & 5.70E-3    & 1.000    & &5.70E-3  &  1.000  &    & 2.38E-6   & 2.000   \\
\hline
\end{tabular}
}
\caption{The $L^2$ errors and rates for Problem 3 with $Y=5,\, \alpha = 0.999999$, $N=30$ and $T=10$.}\label{TabelProblem3}
\end{table}
}}

%%%%%%%%%%%%%%%%%%%%%%%%%%%%%%%%%%%%%%%%%%

{\subsection{Fractional Van der Pol Equation: Problem 4}}
{{As an example of nonlinear problems, we consider the fractional Van der Pol equation \cite{Baffet2017}
%\iffalse%==============>
\begin{equation}
\left(D^{\alpha}\right)^2 x
-\varepsilon (1-x^2)D^{\alpha}x
+x=0,\quad x(0)=x_0,\quad D^\alpha_tx(0)=y_0.
\label{eq:vdp}
\end{equation}
Similar to the previous example, we introduce 
\[y=D^\alpha x, \ x=u+x_0,\ y=v+y_0
\]
to rewrite and discretize  \eqref{eq:vdp} as 
\begin{align}
\label{eq:res-u}
&\delta_tu^{n}+\widehat D_Nu+\sum_{k=1}^N\widehat\sigma_k\phi_{u,k}-\gamma(v+y_0)=0,\\
\label{eq:res-v}
&\delta_tv^{n}+\widehat D_Nv+\sum_{k=1}^N\widehat\sigma_k\phi_{v,k}
    -\gamma\left[\varepsilon\bigl(1-(u+x_0)^2\bigr)(v+y_0)-(u+x_0)\right]=0,\\
\label{eq:res-phiu}
&\delta_t\phi_{u,k}^{n} -\eta_k\phi_{u,k}-u=0,\\
\label{eq:res-phiv}
&\delta_t\phi_{v,k}^{n} -\eta_k\phi_{v,k}-v=0,
\end{align}
where $\gamma=\frac{1}{c_N},\,
    \widehat D_N=\frac{D_N}{c_N},\,
    \widehat\sigma_k=\frac{\sigma_k}{c_N}$ and 
$$
\delta_tz^{n}=\frac{3z^n-4z^{n-1}+z^{n-2}}{2h}.
$$
For the nonlinear problem, the solver used in the previous examples is no longer applicable and the Newton method is applied to address the nonlinearity. More precisely, at the $m$-th iteration, the linearized system reads 
%\fi%<==============
\begin{equation} \label{eq:blockbdf_newton_system}
\begin{bmatrix}
 \mathbf{A}_1 & \mathbf{B}\\
 \mathbf{E}   & \mathbf{C}
\end{bmatrix}
\begin{bmatrix}
\Delta \mathbf{U}^{(m)}\\ 
\Delta \mathbf{W}^{(m)}
\end{bmatrix}
=
-\begin{bmatrix}
 \mathbf{R^{(m)}_U} \\ 
 \mathbf{R^{(m)}_W}
\end{bmatrix}, \quad
\mathbf{U}=\begin{bmatrix}u\\ v\end{bmatrix},\quad
\mathbf{W}=\begin{bmatrix}\phi_u\\ \phi_v\end{bmatrix},
%%%%
\end{equation}
where $\mathbf{B,\,E,\,C}$ are given in the previous example and 
\begin{equation*}
 \mathbf{A}_1=
 \begin{bmatrix}
  a&-\gamma \\ 
  p^{(m)}&q^{(m)}
 \end{bmatrix},\
 p^{(m)}=\frac{1}{c_N}\left(2\varepsilon x^{(m)}y^{(m)}+1\right),
 \
 q^{(m)} = a-\gamma\varepsilon\bigl(1-(x^{(m)})^2\bigr).
\end{equation*}
The Schur complement is
$$
\mathbf{S} =\mathbf{A}_1-\mathbf{B}\mathbf{C}^{-1}\mathbf{E},
$$ 
which leads to 
$$
\Delta\mathbf{U}^{(m)}= \mathbf{S}^{-1}\left(- \mathbf{R_U}+\mathbf{B}\mathbf{C}^{-1}\mathbf{R_W} \right).
$$
}}

{In the implementation, we set $\epsilon=4,\,x_0=2,\,y_0=0,\,T=12,\,\alpha=0.8$ and the tolerance of the Newton method is $10^{-10}$. Because there is no analytical expression for the exact solution of (\ref{eq:vdp}), to study the convergence order, we use the numerical solution $x^r$ corresponding to $M =64000$ and $Y=5,\,N=70$ as the reference solution.}

\begin{figure}[htbp]
\centering
    \begin{subfigure}{0.48\linewidth}
        \centering
        \includegraphics[width=\linewidth]{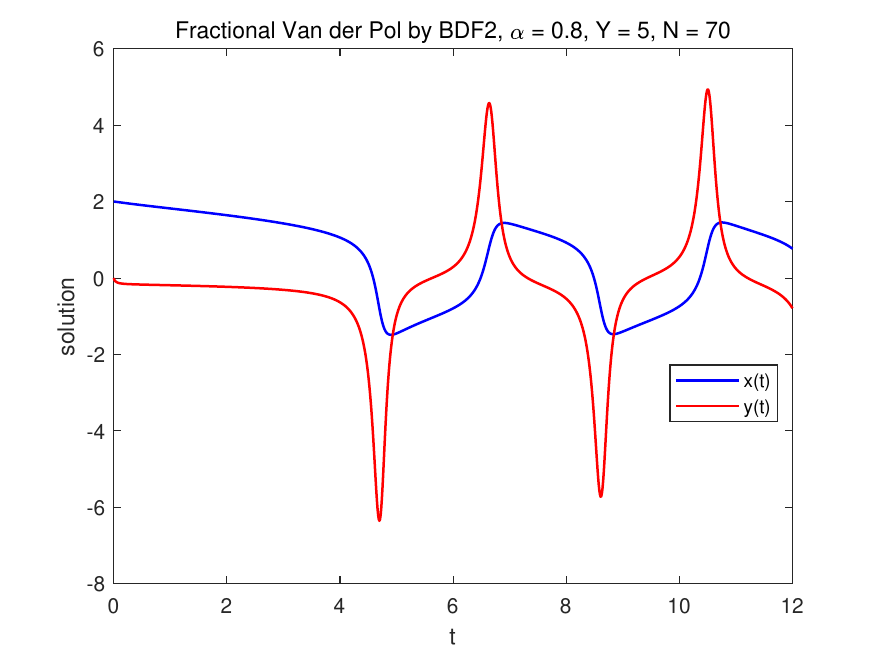}
        \caption{}
    \end{subfigure}\hfill
    \begin{subfigure}{0.48\linewidth}
        \centering
        \includegraphics[width=\linewidth]{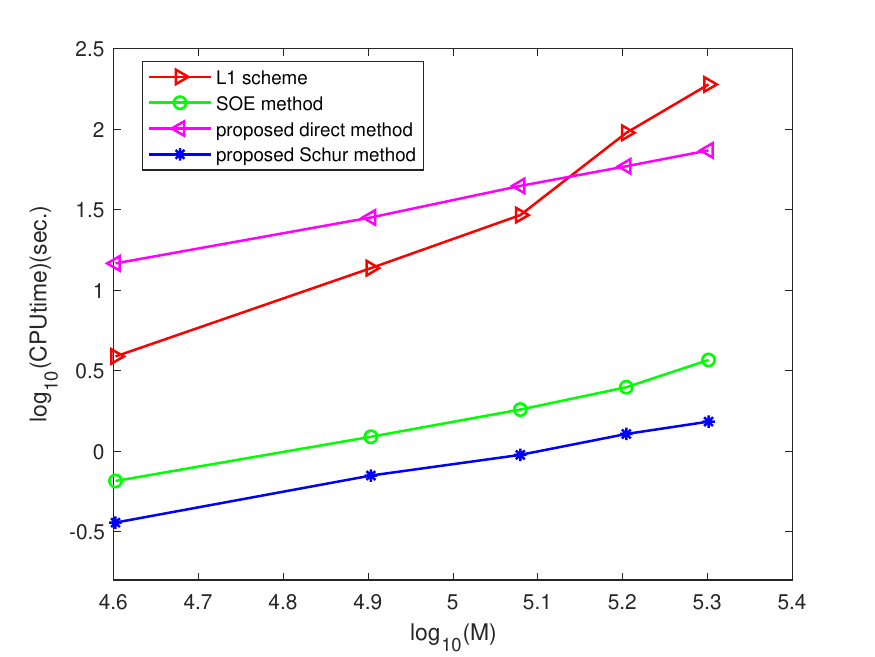}
        \caption{}
    \end{subfigure}
    \caption{Numerical results for the Van der Pol problem. (a) the reference solution with $M=64000$, (b) the CPU time with different $M$, where {\it{proposed direct method }} means solving the problem without the Schur complement. 
    \label{fig:Vanderpol}}
\end{figure}

{The reference solution $x^r$ is shown in Figure \ref{fig:Vanderpol} (a); the result is consistent with the solution in \cite{Baffet2017}.}

{In Table \ref{TabelProblem4}, we report the $L^2$ errors and convergence orders. The expected rate of second order is observed for the proposed method, while the $L_1$ method and the SOE method exhibit a convergence order about $\mathcal{O}(2-\alpha)$.}

{Figure \ref{fig:Vanderpol} (b) shows the CPU time with $M=40000,...,200000$, an expected slope of 1 and a significant reduction in CPU time of the proposed method are observed. We should point out that a direct implementation of the proposed method without the Shur complement, i.e. using \eqref{eq:BDF2} instead of \eqref{eq:BDF2solution}, is more expensive than the SOE method. 
}

\begin{table}
\centering
{
\begin{tabular}{p{10pt} p{5pt} p{50pt} p{20pt}  p{5pt} p{50pt} p{20pt} p{5pt} p{50pt} p{20pt} }
\hline
   \multirow{2}{*}{M}& & \multicolumn{2}{c}{$L_1\, method$}&  & \multicolumn{2}{c}{$SOE\, method$}& & \multicolumn{2}{c}{$Proposed\, method$} \\
\cline{3-4} \cline{6-7} \cline{9-10}
  & & Error &Order     & &    Error &Order &    & Error &Order \\
\hline
 $1000$  &  &1.875E-1    & -         & &1.875E-1  & -       &    & 3.897E-2   & - \\
 $2000$  &  &7.505E-2    & 1.321     & &7.505E-2  &  1.321  &    & 9.953E-3   & 1.969\\
 $4000$  &  &2.902E-2    & 1.371     & &2.902E-2  &  1.371  &    & 2.505E-3   & 1.990 \\
 $8000$  &  &1.073E-3    & 1.436     & &1.073E-3  &  1.436  &    & 6.228E-4   & 2.008 \\
 $16000$ &  &3.665E-3    & 1.550     & &3.665E-3  &  1.550  &    & 1.492E-4   & 2.061   \\
\hline
\end{tabular}
}
\caption{The $L^2$ errors and rates for Problem 4 with $Y=5,\, \alpha = 0.8$, $N=70$ and $T=12$.}\label{TabelProblem4}
\end{table}

%%%%%%%%%%%%%%%%%%%%%%%%%%%%%%%%%%%%%%%%%%
\section{Conclusion {and discussion} \label{Conclusion}}
In this paper, an analysis-based numerical algorithm is proposed to handle the Caputo fractional derivative through the pole-residue approximation of its kernel in the Laplace transformed domain. A deformed contour satisfying simultaneously  the convergence criteria of the Mellin's formula for inverse Laplace transform and the assumptions of Theorems 1 and  5.1 in \cite{gelfgren1978rational} is used for interpolating the transformed kernel function. Because of this, we are able to analyze how { the approximation error in the $z$-domain translate to the error in the time domain}. The analysis shows that the extrapolation error converges to $0$ exponentially as $t$ increases and that the only part of the error growing with time $t$ is the factor $e^{\eta t}$ in \eqref{error_estimate}, so in practice, we can set $\eta=\frac{1}{T}$ if the simulation is  from $t=0$ up to $t=T$. As is demonstrated by the numerical examples, this approximation works well as long as the Laplace transform of the function which the fractional derivative acts on has all its singularities to the left of the contour. 

{For all the test problems, we applied the BDF2 method for discretization and  compare the results with those obtained by the popular $L_1$ method and the SOE method. In all the test problems, the proposed method demonstrates higher convergence order, less CPU time and lower $L^2$ errors  than the $L_1$ method and the SOE method. In Example 3, the proposed method also demonstrate superior energy conservation property as $\alpha$ gets close to 1. For the first two linear test problems, we also applied the the block exponential method to evaluate the exact solution of the approximate system of \eqref{system_eq_1}and \eqref{eq:auxphiu} to study the error as a function of $N$.}

Most of the sum-of-exponential approximation results in the literature assume the coefficient of $u'$ be $0$ in the approximation and as a result, a large number of terms are required to capture the weakly singular behavior of the kernel function near $0$. For example, hundreds of terms are used in \cite{Baffet2017}. Indeed, this was the main motivation of the work in \cite{10.1093/imanum/drac022}, where the authors proposed an  $[N/N]$ rational approximation of $z^{-\alpha}$  constructed by using the well-known AAA algorithm for interpolating $z^{-\alpha}$ on the \textbf{real interval} $[h,T]=[10^{-3},1]$ of the complex $z$ plane. The intuition behind using  $[N/N]$ instead of $[N-1/N]$ is not immediately clear because $z^{-\alpha}$ is singular at $\infty$ and approaches 0 as $z$ goes to $\infty$. Moreover, the constant in the rational approximation in \cite{10.1093/imanum/drac022} is treated as independent from the poles and residues of the rational function. In contrast, the constant $c_N$ in \eqref{approx_time_domain} is a direct result of the $[N-1/N]$ approximation of $Q(s)$ and is a function of $\gamma_k$ and $\eta_k$, $k=1,\dots,N$, not an independent variable. Another crucial difference between our algorithm and the one used in \cite{10.1093/imanum/drac022} is that our interpolation is carried out for $Q(s)$, a Markov function, whose rational approximation theory {such as that behind Theorem \ref{thm_N_conv}} is well established. More importantly, the interpolation grids in our algorithms are on a contour satisfying the Mellin's criteria needed for the integral to converge to the inverse Laplace transform and hence leads to an explicit and rigorous link between the approximation error in the $z$-domain and that in the time domain.

Even though we restrict the focus in the paper on the Caputo derivative, we expect the proposed method to be easily generalized to many other types of time-fractional derivatives. This is because the causality (the current state only depends on the past) built in the definition implies the Laplace transform of the kernel function is analytic in the upper-half plane. If other property  such as the passivity or power dissipation is also built-in, then the transformed kernel will be a Herglotz-Nevanlinna function \cite{Luger2021On-applications}. Note that the assumption of $u(0)=0$ is not essential because we could have kept it i \eqref{temp_Dalpha} with an additional term $-z^{\alpha-1}u(0)$, whose inverse Laplace transform is $-\frac{t^{-\alpha} u(0)}{\Gamma(1-\alpha)}$, and \eqref{system_eq_1} will be take the form:
\[
c_N\frac{du}{dt} =
\left[F(t,u(t))+ \frac{t^{-\alpha} u_0}{\Gamma(1-\alpha)}\right] - D_N u(t) - \displaystyle\sum_{k=1}^N \sigma_k \phi_k(t),
\qquad u(0)=u_0
\]

{We also remark that the proposed scheme is also applicable for the Caputo derivatives of higher orders $D^{n+\alpha}$, $n\in \mathbb{N}$ and $0<\alpha<1$, defined as  
\[
D^{n+\alpha}u(t):=\frac{1}{\Gamma(1-\alpha)}\int_0^s \frac{u^{(n+1)}(s)}{(t-s)^\alpha} ds
\]
Suppose $n=1$ and the initial conditions are $u(0)=a$ and $u'(0)=b$. Defining $w(t)=u(t)-a-bt$, then $D^{1+\alpha}w=D^{1+\alpha}u$ and $w(0)=0$, $w'(0)=0$. Hence WLOG, we may assume $u(0)=0$ and $u'(0)=0$, which in turn leads to
\[
\mathcal{L}[D^{1+\alpha}u]=z^{\alpha-1}\cdot z^2 \hat{u}
\approx c_N (z^2 \hat u)+\left(D_N+\sum_{k=1}^N \frac{\sigma_k}{z-\eta_k}\right)(z\hat{u}).
\]
Applying the inverse Laplace transform results in
\[
D^{1+\alpha}u(t)\approx c_N u''(t)+D_N u'(t)+\sum_{k=1}^N \sigma_k \int_0^t e^{\eta_k (t-\tau)} u'(\tau) d\tau 
\]
Letting
\beq
\phi_k(t):=\int_0^t e^{\eta_k (t-\tau)} u'(\tau) d\tau,\, k=1,\dots,N
\label{eq_general}
\eeq
then \eqref{eq_general} is equivalent to the following system:
\beqs
&&D^{1+\alpha}u(t)\approx c_N u''(t)+D_N u'(t)+\sum_{k=1}^N \sigma_k \phi_k(t)\\
&&\frac{d\phi_k}{dt}=\sigma_k u'(t)+\eta_k \phi_k,\, \phi_k(0)=0,\,k=1,\dots,N.\\
\eeqs
Generalization to higher $n\in\mathbb{N}$ is straightforward.
}

{
For future work, we will improve the current scheme by addressing the following. First, we will eliminate the need for an arbitrary precision arithmetic system for solving the square linear system of size $2N$ in \eqref{complex_system} by using the continued fraction approximation of Stieltjes functions \cite{celis2024numerical-conti} from which the poles and residues of $Q(s)$ can be extracted directly without first solving \eqref{complex_system}. The other direction we would like to explore is to find a better balance between the two terms in the brackets of \eqref{error_estimate}, the interpolation error and the truncation error, by properly choosing $Y$ as a function of $N$. 
}
\section*{Acknowledgments}
M-J Y. Ou gratefully acknowledges the support of the EU RISE Grant: Project EXPOWER for her work. The work of J. Xie is partially supported by the National Natural Science Foundation of China under Grant 12301550 and Natural Science Foundation of Chongqing under Grant CSTB2024NSCQ-MSX1288.

%%%%%%%%%%%%Citation%%%%%%%%%%%%%%%%%%%%%%%%%%%%%%%%
\bibliographystyle{plain}
\bibliography{library}

@article{10.1093/imanum/drac022,
	author = {Khristenko, U. and Wohlmuth, B.},
	doi = {10.1093/imanum/drac022},
	eprint = {https://academic.oup.com/imajna/article-pdf/43/3/1263/50504270/drac022.pdf},
	issn = {0272-4979},
	journal = {IMA Journal of Numerical Analysis},
	month = {06},
	number = {3},
	pages = {1263-1290},
	title = {Solving time-fractional differential equations via rational approximation},
	url = {https://doi.org/10.1093/imanum/drac022},
	volume = {43},
	year = {2022}}

@book{apostol1974mathematical-an,
	author = {T. Apostol},
	edition = {{S}econd},
	publisher = {Pearson},
	title = {Mathematical Analysis},
	year = {1974}}

@article{dingfelder2015improved,
  title={An improved {T}albot method for numerical {L}aplace transform inversion},
  author={Dingfelder, B. and Weideman, J.},
  journal={Numerical Algorithms},
  volume={68},
  number={1},
  pages={167--183},
  year={2015},
  publisher={Springer}
}

@book{jin2023numerical,
  title={Numerical treatment and analysis of time-fractional evolution equations},
  author={Jin, Bangti and Zhou, Zhi},
  volume={214},
  year={2023},
  publisher={Springer}
}

@inbook{Luger2021On-applications,
	author = {Luger, Annemarie and Ou, Miao-Jung Yvonne},
	pages = {433-459},
	publisher = {Springer},
	series = {AWM Mathematics: Research in the Mathematics of Materials Science},
	title = {On applications of Herglotz-Nevanlinna functions in material sciences, I: classical theory and applications of sum rules},
	volume = {31},
	year = {2022}}

@article{Baffet2017,
	author = {Baffet, D. and Hesthaven, J.},
	doi = {10.1137/15M1043960},
	eprint = {https://doi.org/10.1137/15M1043960},
	journal = {SIAM Journal on Numerical Analysis},
	number = {2},
	pages = {496-520},
	title = {A Kernel Compression Scheme for Fractional Differential Equations},
	url = {https://doi.org/10.1137/15M1043960},
	volume = {55},
	year = {2017}}

@article{celis2024numerical-conti,
	author = {Celis, Oliver Salazar},
	journal = {Ukrainian Mathematical Journal},
	number = {4},
	pages = {635--648},
	title = {Numerical Continued Fraction Interpolation},
	volume = {76},
	year = {2024}}

@article{article,
	author = {Aghili, Arman and Ansari, Alireza},
	journal = {Int. J. Contemp. Math. Sciences},
	month = {01},
	pages = {761-770},
	title = {Complex inversion formula for {S}tieltjes and {W}idder transforms with applications},
	volume = {3},
	year = {2008}}

@article{Daniel-Alpay1994The-two-sided-r,
	author = {Alpay, D. and Ball, J. and Gohberg, I. and Rodman, L.},
	journal = {Linear Algebra and Its Applications},
	pages = {485-521},
	title = {The two-sided residue interpolation in the {S}tieltjes class for matrix functions},
	volume = {208-209},
	year = {1994}}

@article{Kac1974R-functions-ana,
	author = {Kac, I. and Krein, M.},
	journal = {AMS Translations},
	pages = {1-18},
	title = {R-functions-analytic functions mapping the upper halfplane into itself},
	volume = {103},
	year = {1974}}

@article{VKStieltjesFunctionsandHurwitz,
	author = {Katsnelson, V.},
	doi = {10.1007/s11785-011-0146-1},
	journal = {Complex Analysis and Operator Theory},
	title = {{S}tieltjes functions and {H}urwitz stable entire functions},
	volume = {5},
         pages={611--630},
	year = {2011}}

@incollection{ou2018augmented,
	author = {Ou, M. and Woerdeman, H.},
	booktitle = {Interpolation and Realization Theory with Applications to Control Theory In Honor of Joe Ball},
	editor = {Vladimir Bolotnikov and Sanne ter Horst and Andr\'{e} C.M. Ran and Victor Vinnikov},
	journal = {Operator Theory: Advances and Applications, Springer Nature},
	pages = {307-328},
	publisher = {Springer},
	series = {Operator Theory: Advances and Applications},
	title = {On the augmented {B}iot-{JKD} equations with Pole-Residue representation of the dynamic tortuosity},
	volume = {272},
	year = {2019}}

@article{ou2014on-reconstructi,
	author = {Ou, M.},
	journal = {Inverse Problems},
	number = {9},
	pages = {095002},
	title = {On reconstruction of dynamic permeability and tortuosity from data at distinct frequencies},
	volume = {30},
	year = {2014}}

@techreport{gelfgren1978rational,
	author = {Gelfgren, J.},
	institution = {Department of Mathematics, University of Ume{\aa}},
	title = {Rational interpolation to functions of {S}tieltjes' type},
	year = {1978}}

@incollection{gelfgren1981multipoint,
	author = {Gelfgren, J.},
	booktitle = {Pad{\'e} Approximation and its Applications Amsterdam},
	pages = {197--207},
	publisher = {Springer},
	title = {Multipoint {P}ad{\'e} approximants converging to functions of {S}tieltjes' type},
	year = {1981}}

@book{akhiezer1993theory-of-linea,
	author = {Akhiezer, N. and Glazman, I.},
	publisher = {Dover Publications},
	title = {Theory of linear operators in Hilbert space},
	volume = {1},
	year = {1993}}

@article{IdizTanogluAghazadeh2024,
  author  = {Idiz, F. and Tanoglu, G. and Aghazadeh, N.},
  title   = {A numerical method based on {L}egendre wavelet and quasilinearization technique for fractional {L}ane--{E}mden type equations},
  journal = {Numerical Algorithms},
  year    = {2024},
    volume = {95},
   pages = {181–206}}

@article{AghazadehMohammadiTanoglu2024,
  author  = {Aghazadeh, N. and Mohammadi, A. and Tanoglu, G.},
  journal = {Mathematical Sciences},
	pages = {41–54},
   title   = {Taylor wavelets collocation technique for solving fractional nonlinear singular {P}{D}{E}s},
   volume = {18},
     year    = {2024},
}

@article{2010Fractional,
  title={Fractional Calculus and Waves in Linear Viscoelasticity: An Introduction to Mathematical Models},
  author={ Mainardi, F. },
  journal={Imperial College Press},
  year={2010},
}

@book{kilbas2006theory,
  title={Theory and applications of fractional differential equations},
  author={Kilbas, A. and Srivastava, H. and Trujillo, J.},
  volume={204},
  year={2006},
  publisher={{E}lsevier}
}

@article{metzler2000random,
  title={The random walk's guide to anomalous diffusion: a fractional dynamics approach},
  author={Metzler, R. and Klafter, J.},
  journal={Physics reports},
  volume={339},
  number={1},
  pages={1--77},
  year={2000},
  publisher={Elsevier}
}

@book{samko1993fractional,
  title={Fractional integrals and derivatives},
  author={Samko, S. and A. Kilbas and O. I. Marichev},
  year={1993},
  publisher={Gordon and Breach}
}

@article{diethelm2002analysis,
  title={Analysis of fractional differential equations},
  author={Diethelm, K. and Ford, N.},
  journal={Journal of Mathematical Analysis and Applications},
  volume={265},
  number={2},
  pages={229--248},
  year={2002},
  publisher={Elsevier}
}

@book{podlubny1998fractional,
  title={Fractional differential equations: an introduction to fractional derivatives, fractional differential equations, to methods of their solution and some of their applications},
  author={Podlubny, I.},
  volume={198},
  year={1998},
  publisher={{E}lsevier}
}

@article{lubich1988convolution,
  title={Convolution quadrature and discretized operational calculus. {I}},
  author={Lubich, C.},
  journal={Numerische Mathematik},
  volume={52},
  number={2},
  pages={129--145},
  year={1988},
  publisher={Springer}
}

@article{jin2018discrete,
  title={Discrete maximal regularity of time-stepping schemes for fractional evolution equations},
  author={Jin, B. and Li, B. and Zhou, Z.},
  journal={Numerische mathematik},
  volume={138},
  number={1},
  pages={101--131},
  year={2018},
  publisher={Springer}
}

@article{jiang2017fast,
  title={Fast evaluation of the {C}aputo fractional derivative and its applications to fractional diffusion equations},
  author={Jiang, S. and Zhang, J. and Zhang, Q. and Zhang, Z.},
  journal={Communications in Computational Physics},
  volume={21},
  number={3},
  pages={650--678},
  year={2017},
  publisher={Cambridge University Press}
}

@article{banjai2019efficient,
  title={Efficient high order algorithms for fractional integrals and fractional differential equations},
  author={Banjai, L. and L{\'o}pez-Fern{\'a}ndez, M.},
  journal={Numerische Mathematik},
  volume={141},
  number={2},
  pages={289--317},
  year={2019},
  publisher={Springer}
}

@article{schadle2006fast,
  title={Fast and oblivious convolution quadrature},
  author={Sch{\"a}dle, A. and L{\'o}pez-Fern{\'a}ndez, M. and Lubich, C.},
  journal={SIAM Journal on Scientific Computing},
  volume={28},
  number={2},
  pages={421--438},
  year={2006},
  publisher={SIAM}
}

@article{yuan2002numerical,
  title={A numerical scheme for dynamic systems containing fractional derivatives},
  author={Yuan, L. and Agrawal, O.},
  journal={Journal of Vibration and Acoustics},
  volume={124},
  number={2},
  pages={321--324},
  year={2002}
}

@article{lubich2002fast,
  title={Fast convolution for nonreflecting boundary conditions},
  author={Lubich, C. and Sch{\"a}dle, A.},
  journal={SIAM Journal on Scientific Computing},
  volume={24},
  number={1},
  pages={161--182},
  year={2002},
  publisher={SIAM}
}

@article{rosenfeld2017approximating,
  title={Approximating the {C}aputo fractional derivative through the {M}ittag-{L}effler reproducing kernel {H}ilbert space and the kernelized {A}dams--{B}ashforth--{M}oulton method},
  author={Rosenfeld, J. and Dixon, W.},
  journal={SIAM Journal on Numerical Analysis},
  volume={55},
  number={3},
  pages={1201--1217},
  year={2017},
  publisher={SIAM}
}

@article{yuste2005explicit,
  title={An explicit finite difference method and a new von {N}eumann-type stability analysis for fractional diffusion equations},
  author={Yuste, S. and Acedo, L.},
  journal={SIAM Journal on Numerical Analysis},
  volume={42},
  number={5},
  pages={1862--1874},
  year={2005},
  publisher={SIAM}
}

@article{zayernouri2016fractional,
  title={Fractional {A}dams--{B}ashforth/{M}oulton methods: an application to the fractional {K}eller--{S}egel chemotaxis system},
  author={Zayernouri, M. and Matzavinos, A.},
  journal={Journal of Computational Physics},
  volume={317},
  pages={1--14},
  year={2016},
  publisher={Elsevier}
}

@book{oldham1974fractional,
  title={The fractional calculus theory and applications of differentiation and integration to arbitrary order},
  author={Oldham, K. and Spanier, J.},
  volume = {111},
  year={1974},
  publisher={Elsevier}
}

@article{xie2019discontinuous,
  title={A discontinuous {G}alerkin method for wave propagation in orthotropic poroelastic media with memory terms},
  author={Xie, J. and Ou, M. and Xu, L.},
  journal={Journal of Computational Physics},
  volume={397},
  pages={108865},
  year={2019},
  publisher={Elsevier}
}

@article{zeng2015numerical,
  title={Numerical algorithms for time-fractional subdiffusion equation with second-order accuracy},
  author={Zeng, F. and Li, C. and Liu, F. and Turner, I.},
  journal={SIAM Journal on Scientific Computing},
  volume={37},
  number={1},
  pages={A55--A78},
  year={2015},
  publisher={SIAM}
}

@article{zheng2023fast,
  title={Fast artificial boundary method for the heat equation on unbounded domains with strip tails},
  author={Zheng, C. and Xie, J.},
  journal={Journal of Computational and Applied Mathematics},
  volume={425},
  pages={115032},
  year={2023},
  publisher={Elsevier}
}
%%%%%%

\end{document}